\documentclass{article} 
\usepackage{paper2lean_arxiv,times}

\usepackage{amsmath,amsfonts,bm}

\def\eqref#1{equation~\ref{#1}}

\def\1{\bm{1}}

\DeclareMathAlphabet{\mathsfit}{\encodingdefault}{\sfdefault}{m}{sl}
\SetMathAlphabet{\mathsfit}{bold}{\encodingdefault}{\sfdefault}{bx}{n}

\usepackage{hyperref}
\usepackage{url}
\usepackage{multirow}
\usepackage{listings}

\usepackage{graphicx}
\usepackage{booktabs}
\usepackage{tikz}
\usepackage{listings}
\usepackage{tcolorbox}
\usetikzlibrary{arrows.meta,positioning,fit,calc}
\title{Paper2Lean: Agentic Autoformalization of Mathematical Research Papers at Scale}

\iclrfinalcopy

\author{Yifan Bai\thanks{These authors contributed equally.},
Yantao Li\footnotemark[1],
Jian Yu\footnotemark[1] \\
School of Mathematical Sciences \\
Shanghai Jiao Tong University, China \\
\texttt{\{yifanbai,alexandero,yujian26\}@sjtu.edu.cn} \\
\And
Jingwei Liang\thanks{Corresponding author.} \\
Institute of Natural Sciences \\
Shanghai Jiao Tong University, China \\
\texttt{jingwei.liang@sjtu.edu.cn}
}

\usepackage{xcolor}
\usepackage{listings}

  \definecolor{jsonBackground}{HTML}{F7F9FC}
  \definecolor{jsonFrame}{HTML}{CBD5E1}

  \lstdefinestyle{jsonArtifact}{
    basicstyle=\ttfamily\scriptsize\color{black},
    stringstyle=\color{black},
    keywordstyle=\color{black},
    morestring=[b]",
    morekeywords={true,false,null},
    backgroundcolor=\color{jsonBackground},
    rulecolor=\color{jsonFrame},
    frame=single,
    framerule=0.5pt,
    framesep=6pt,
    xleftmargin=0.025\linewidth,
    xrightmargin=0.025\linewidth,
    framexleftmargin=5pt,
    framexrightmargin=5pt,
    framextopmargin=4pt,
    framexbottommargin=4pt,
    breaklines=true,
    breakatwhitespace=false,
    breakindent=1.2em,
    columns=fullflexible,
    keepspaces=true,
    showstringspaces=false,
    captionpos=b,
    aboveskip=11pt plus 2pt minus 2pt,
    belowskip=10pt plus 2pt minus 2pt,
    abovecaptionskip=7pt,
    belowcaptionskip=0pt
  }

\begin{document}

\maketitle

\begin{abstract}
Formalizing a mathematical research paper requires more than translating isolated theorems: statements and proofs must remain semantically faithful to the source and globally consistent across long-range dependencies. We introduce \textsc{Paper2Lean}, an end-to-end agentic framework that converts mathematical research papers into integrated Lean projects. The framework extracts mathematical content into a dependency graph and coordinates specialized agents for statement formalization, semantic checking and repair, proof construction, and project integration. We evaluate \textsc{Paper2Lean} on 20 mathematical optimization research papers, using \texttt{mathlib} and a shared library of supporting results. All 20 automatically generated projects compile successfully without proof placeholders. Following human review and semantic repair, we introduce \textsc{OptPaper-20}, a corpus of 20 corrected paper-level formalizations comprising 494 Lean files and 110,910 lines of Lean code. All released projects pass clean builds, with no unresolved discrepancies identified in the audit. These results demonstrate the feasibility of agentic paper-level formalization while highlighting semantic fidelity as a central challenge beyond proof verification. The source code and dataset are publicly available at \url{https://github.com/TTony2019/OptPaper}.
\end{abstract}

\section{Introduction}

Interactive theorem provers such as Rocq~\citep{barras1997coq}, Isabelle~\citep{nipkow2002isabelle}, and
Lean~\citep{de2015lean} have become essential infrastructure for AI for Mathematics
(AI4Math), providing formal languages and trusted kernels for constructing
and verifying machine-checkable mathematics. Among them, Lean is based on dependent type theory and supports the
representation of mathematical statements and proofs as executable,
machine-checkable code. Its mathematical library, \texttt{mathlib}
\citep{mathlib2020}, offers an extensive and continuously growing collection
of reusable definitions, theorems, and proof automation across diverse areas
of mathematics. Together, Lean and \texttt{mathlib} provide a particularly suitable
foundation for research-level formalization: once a mathematical result has
been faithfully encoded, its proof can be checked by Lean's small trusted
kernel and subsequently reused in larger formal developments.


AI for mathematics has advanced rapidly in recent years. In 2026, AI-assisted systems made particularly notable progress 
in conjecture resolution and formal proof search, including mathematical
results checked by Lean~\citep{ju2026automated,tsoukalas2026advancing}.
Large scale efforts also released complete Lean formalizations of major
mathematical results, such as Fermat's Last Theorem~\citep{anthropic2026flt}.
These developments show that LLM agents can contribute not only to discovering or
proving individual results, but also to producing reusable, machine-checked
mathematical knowledge. The next challenge is to formalize the research papers
that contain such results as coherent Lean projects.

Formalizing an entire research paper is substantially different from
translating an isolated theorem. First, the
mathematical content must be recovered from a PDF, including notation,
cross-references, and the structure of definitions and results. Second, the
formal statements must remain aligned with the source despite omitted or
implicit assumptions and variations in notation. Third, the resulting proofs
and declarations must be integrated into one project whose imports,
interfaces, and downstream dependencies remain consistent as the formalization
evolves.

\begin{figure}[ht]
\centering\includegraphics[width=1\linewidth]{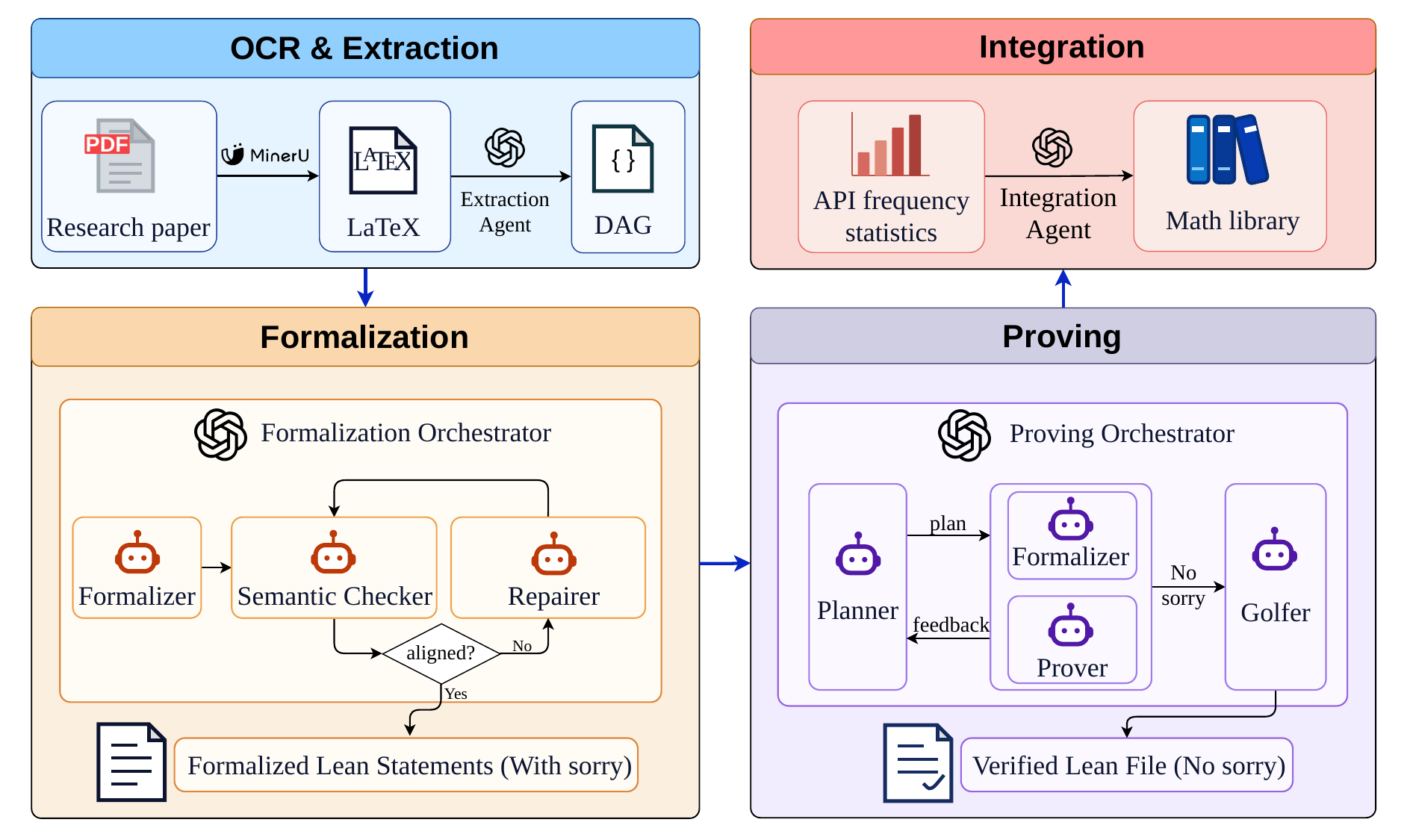}
    \caption{Overview of \textsc{Paper2Lean}.}
    \label{fig:pipeline}
\end{figure}

To address these challenges, we introduce \textsc{Paper2Lean}, an end-to-end agentic framework for formalizing research papers as integrated Lean projects. Figure~\ref{fig:pipeline} provides an overview.
\textsc{Paper2Lean} first applies OCR and
document extraction to identify the mathematical content of a paper and
organize it into a structured dependency graph. It then coordinates specialized
agents for statement formalization, semantic checking, library API reuse, proof
planning and construction, proof refactoring, and library integration.
Compilation gates are applied throughout the pipeline, and verified interfaces
are frozen before dependent items are processed. This design allows independent
items to be handled concurrently while maintaining consistent definitions,
statements, proofs, and dependencies across the complete paper.

We evaluate \textsc{Paper2Lean} on 20 research papers in mathematical
optimization. The automatically generated projects all compile successfully and
contain no \texttt{sorry} placeholders. Experts
judged 473 of 494 items to be source-faithful (95.7\%). After manual correction, we
release the corrected project as \textsc{OptPaper-20}, a corpus containing
494 Lean files and 110,910 lines of Lean code.
All identified discrepancies were repaired in the released projects.

Our contributions are summarized as follows:
\begin{enumerate}
\item We introduce \textsc{Paper2Lean}, an end-to-end agentic framework for
formalizing mathematical research papers in Lean. Starting from a paper PDF,
\textsc{Paper2Lean} extracts and organizes mathematical content into a
dependency graph, then coordinates specialized agents for statement
formalization, semantic checking, proof construction,
golfing, and project integration. Compilation-gated execution and frozen
interfaces enable parallel processing while preserving semantic fidelity,
Lean-checked correctness, and global consistency across the paper.
\item We evaluate \textsc{Paper2Lean} on 20 research papers in mathematical
optimization. All generated projects compile without proof placeholders. An
expert audit finds that 473 of 494 formalized items are source-faithful, while the remaining discrepancies demonstrate that Lean
compilation alone does not establish fidelity to the source paper.
\item We release \textsc{OptPaper-20}, an expert-refined dataset built from the
formalizations of these 20 papers. For each paper, the dataset provides the
refined Lean files, corresponding natural language descriptions, and the
structured JSON file extracted from the source paper, together with mappings
between source paper items and their formal counterparts.
\end{enumerate}

\section{Related Work}

\subsection{Autoformalization and automated theorem proving}

LLM-based autoformalization can be broadly divided into statement-level and
proof-level formalization. Statement-level methods translate natural language
mathematical statements into Lean declarations. They use tools such as \texttt{mathlib}
retrieval and Lean feedback to improve formalizations and construct high quality
datasets for fine tuning, see \citep{hanyu2026aria, liu2026decompose,liu2026atlas}.
Proof-level methods formalize complete
natural language proofs at the project level. ProofFlow~\citep{cabral2026proofflow}
organizes dependencies among informal proof steps as a directed acyclic graph (DAG) 
and formalizes the resulting intermediate lemmas and proofs into a Lean project.
Our work instead focuses on formalizing an entire mathematical paper as a single Lean project, guided by the DAG of dependencies among the paper's theorems.


Automated theorem proving assumes a formally specified theorem and aims to
construct a machine-checkable proof. Recent neural theorem provers combine
LLM generation with premise retrieval, proof search, and Lean feedback, such as LeanDojo and ReProver~\citep{yang2023leandojo}. More recent systems, such as
DeepSeek-Prover-V2~\citep{ren2025deepseek} and Kimina-Prover~\citep{wang2025kimina}, improve proof generation through
subgoal decomposition, structured reasoning, and reinforcement learning.
While these methods achieve strong performance on isolated theorem-proving
benchmarks~\citep{zheng2021minif2f,azerbayev2023proofnet}, they generally assume that theorem statements and their surrounding
formal contexts have already been provided, rather than constructing a
consistent formal development from a mathematical document.

\subsection{Agentic formal reasoning}

Recent work has explored agentic workflows that coordinate large language models,
formal verifiers, retrieval systems, and specialized reasoning modules.
Numina-Lean-Agent~\citep{liu2026numina} equips a coding agent with Lean
interaction, theorem retrieval, and iterative proof-repair tools, while
M2F~\citep{wang2026m2f} extends agentic formalization to entire documents and
projects through dependency-aware statement construction and proof repair.
For research-level mathematics, Archon~\citep{ju2026automated} decomposes
informal arguments and translates them into multi-file Lean projects using
retrieval and verifier-guided refinement. 
Together, these systems demonstrate the benefits of tool use, task
decomposition, and specialized-agent collaboration. Our work focuses on paper-level formalization, where agents must maintain
global consistency and resolve long-range dependencies among definitions,
theorems, and proofs, rather than solving each theorem in isolation.

\section{Paper2Lean}
\label{sec:pipeline}

In this section, we describe the main stages of \textsc{Paper2Lean} shown in Figure~\ref{fig:pipeline}, including OCR and extraction, statement formalization, theorem proving and library integration.

\subsection{OCR and extraction}


Given a research paper in PDF form, we use MinerU~\citep{wang2024mineru} to obtain Markdown and LaTeX versions of the document.
The
extraction agent then reads the entire document and identifies the mathematical units required for formalization, including standing assumptions, definitions, algorithms, lemmas, propositions, theorems, corollaries, and formally relevant remarks. Results cited from external sources are represented
explicitly. The
agent also introduces definition-gap items when a construction is used by the paper but not explicitly defined, and propagates section-level assumptions to all affected results. 
The extracted content is stored in a canonical JSON representation. Each item corresponds to one mathematical unit and records its source statement, optional proof text, type, and dependencies. See Appendix~\ref{appendix:intermediate-artifact} for an example.

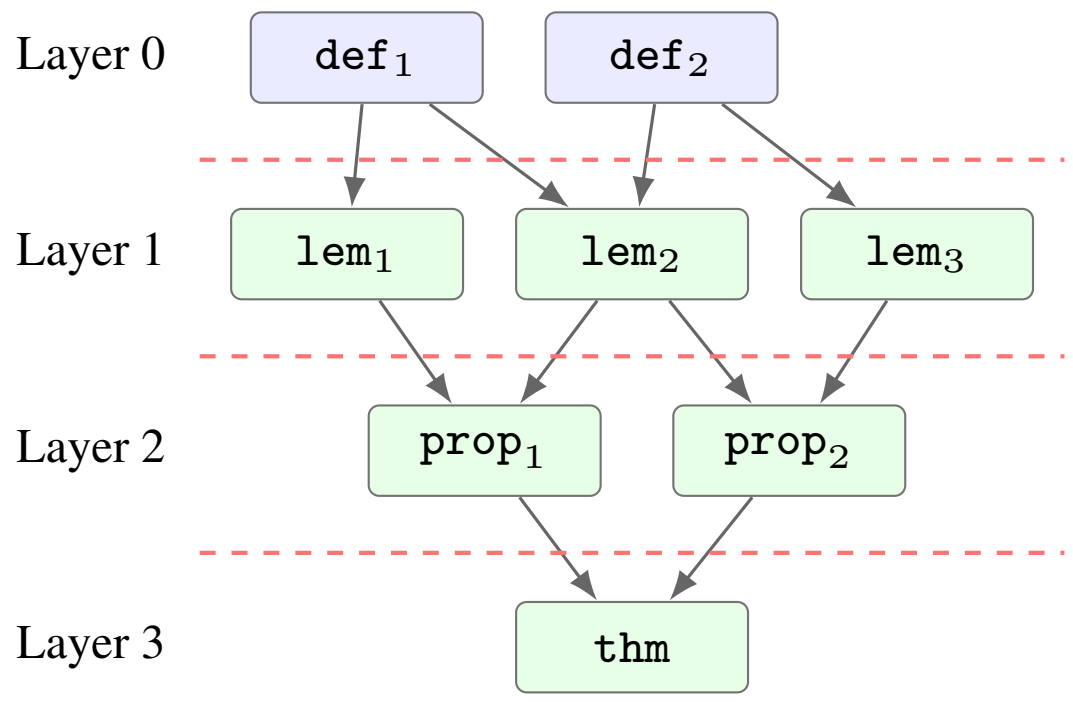
\begin{figure}[h]
\centering
\begin{tikzpicture}[font=\scriptsize,>=Latex,scale=1.3,transform shape,
  item/.style={draw=black!55,fill=green!10,rounded corners=2pt,minimum width=1.18cm,minimum height=0.46cm,align=center},
  dep/.style={draw=black!55,fill=blue!8,rounded corners=2pt,minimum width=1.18cm,minimum height=0.46cm,align=center},
  barrier/.style={draw=red!55,dashed,thick},
  arr/.style={->,semithick,black!60}]
  \node[dep] (d1) at (0,0) {${\tt def}_1$};
  \node[dep] (d2) at (1.5,0) {${\tt def}_2$};
  \node[item] (l1a) at (-0.1,-1.0) {${\tt lem}_1$};
  \node[item] (l1b) at (1.35,-1.0) {${\tt lem}_2$};
  \node[item] (l1c) at (2.8,-1.0) {${\tt lem}_3$};
  \node[item] (l2a) at (0.6,-2.0) {${\tt prop}_1$};
  \node[item] (l2b) at (2.15,-2.0) {${\tt prop}_2$};
  \node[item] (l3) at (1.35,-3.0) {${\tt thm}$};
  \draw[arr] (d1) -- (l1a);
  \draw[arr] (d1) -- (l1b);
  \draw[arr] (d2) -- (l1b);
  \draw[arr] (d2) -- (l1c);
  \draw[arr] (l1a) -- (l2a);
  \draw[arr] (l1b) -- (l2a);
  \draw[arr] (l1b) -- (l2b);
  \draw[arr] (l1c) -- (l2b);
  \draw[arr] (l2a) -- (l3);
  \draw[arr] (l2b) -- (l3);
  \draw[barrier] (-0.85,-0.52) -- (3.55,-0.52);
  \draw[barrier] (-0.85,-1.52) -- (3.55,-1.52);
  \draw[barrier] (-0.85,-2.52) -- (3.55,-2.52);
  \node[anchor=east] at (-0.9,0) {Layer 0};
  \node[anchor=east] at (-0.9,-1) {Layer 1};
  \node[anchor=east] at (-0.9,-2) {Layer 2};
  \node[anchor=east] at (-0.9,-3) {Layer 3};
\end{tikzpicture}
\caption{Illustrative DAG for paper-level formalization. The labels
indicate item types with {\tt def} for definition, {\tt lem} for lemma, {\tt prop} for proposition,
and {\tt thm} for theorem. 
A lemma may depend on a definition, while a theorem may depend
on preceding propositions; dashed lines indicate the layer barriers used by the
proving scheduler.
}
\label{fig:dag-scheduling}
\end{figure}

Following the approach proposed by \citet{ju2026matlas}, we represent each paper as a directed acyclic graph (DAG), where each node corresponds to a mathematical unit and directed edges encode dependencies between units.
A directed edge from $u$ to $v$ means that the declaration or proof of $v$ is based on
$u$ (see Figure~\ref{fig:dag-scheduling}); consequently, prerequisites appear
above the items that depend on them. The
upper layers therefore typically contain standing assumptions and basic
definitions, while lower layers contain lemmas and theorems that build on them.
The extraction and planning stages validate the references and reject
cycles, so every item can be assigned a dependency depth.

\subsection{Statement formalization}

The formalization orchestrator coordinates three workers for each item: the \textbf{formalizer}, the \textbf{semantic checker}, and the \textbf{semantic repairer}.
For each item, the formalizer first reads the source statement and its dependencies. It then imports the corresponding dependency modules, searches \texttt{mathlib} and the shared \texttt{Foundation} for relevant APIs, avoiding redundant definitions and improving compatibility with the existing formalization ecosystem. An initial Lean formalization of the original statement is then generated.

Since the generated statement may not fully preserve the intended meaning of its natural language counterpart, it is subsequently reviewed by the semantic checker. If a semantic misalignment is detected, the checker provides detailed feedback, marks the item with \texttt{inaccurate}, and the orchestrator assigns the statement to the semantic repairer to fix the issue. The repaired statement is then returned to the semantic checker for another round of validation. This iterative repair-and-check loop continues until the checker confirms that the formal statement is semantically aligned with the source statement. Statements are frozen before the proving stage unless the proving agent identifies a gap in the statement or its assumptions.

\subsection{Theorem proving}

The theorem proving stage replaces the remaining \texttt{sorry} placeholders with Lean proofs. The \textbf{proving orchestrator} coordinates four roles: \textbf{planner} (proof planning),  \textbf{formalizer} (statement formalization),  \textbf{prover} (proof generation) and \textbf{golfer} (proof golfing).  The same DAG is used again, 
items in one depth layer can be proved
concurrently, whereas a dependent layer waits until every item in the current layer
has a compiling, frozen interface. This barrier preserves a stable context for
parallel proving agents to make
progress simultaneously.


For each item, the planner first analyzes the target \texttt{sorry} placeholders and reads the original natural language proof. It then develops a proof plan by decomposing the target into manageable subgoals, identifying potentially useful lemmas and tactics, and determining an appropriate proof strategy.
Based on this plan, the formalizer states useful auxiliary lemmas and leaves their proofs as placeholders. These lemmas are kept in the same Lean file as the target theorem and are not added to the paper-level DAG.
The prover is asked to construct the corresponding Lean proof and fill in the \texttt{sorry}. If the proof attempt fails, 
the prover returns diagnostic information and structured feedback to the orchestrator. The orchestrator then uses this feedback to request a revised proof plan and ask for another proof attempt. This plan–prove–feedback loop continues until all proof obligations have been discharged and the resulting Lean file contains no \texttt{sorry}. Only after repeated unsuccessful plan–prove iterations, when the planner explicitly identifies missing assumptions, may the proving orchestrator add the corresponding assumptions to the statement. An item is considered proved only when it contains no \texttt{sorry}, \texttt{admit}, \texttt{axiom} or Lean errors.

Finally, the completed proof is passed to the subsequent proof-golfing stage, where it can be simplified and improved without changing its correctness. Our golf stage is based on Lean skills~\footnote{\url{https://github.com/cameronfreer/lean4-skills}}. It also runs the
\texttt{mathlib} linters and revises the code according to their
diagnostics, including warnings about redundant arguments, deprecated
declarations, unnecessary tactics, and noncanonical style.

\subsection{Library integration}

Library integration is performed after a batch of papers has been formalized.
The integration agent first groups candidate
definitions and theorems expressing the same mathematical content, and
summarizes their occurrence and use across distinct papers. Frequently
recurring candidates are prioritized for consolidation into shared interfaces, keeping \texttt{Foundation} focused on reusable content
and limiting the accumulation of redundant or narrowly specialized declarations.

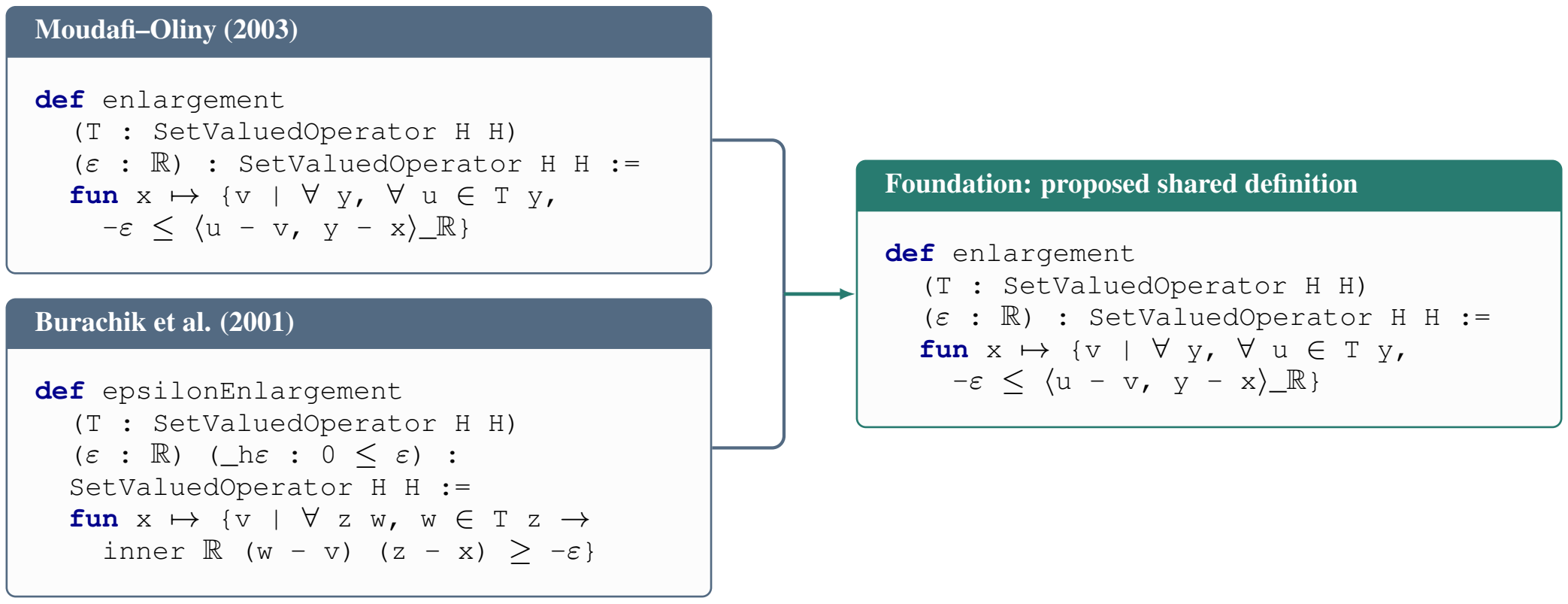
\begin{figure}[htbp]
\centering
\begingroup
\definecolor{enlargeLocal}{HTML}{536A82}
\definecolor{enlargeShared}{HTML}{287B70}
\lstdefinestyle{enlargementLean}{
  basicstyle=\ttfamily\scriptsize,
  keywordstyle=\color{blue!55!black}\bfseries,
  morekeywords={def,fun},
  columns=fullflexible,keepspaces=true,
  showstringspaces=false,breaklines=true,
  aboveskip=0pt,belowskip=0pt,
  literate={ℝ}{{$\mathbb{R}$}}1 {ε}{{$\varepsilon$}}1
           {∀}{{$\forall$}}1 {∈}{{$\in$}}1
           {≤}{{$\leq$}}1 {≥}{{$\geq$}}1
           {→}{{$\rightarrow$}}1 {↦}{{$\mapsto$}}1
           {⟪}{{$\langle$}}1 {⟫}{{$\rangle$}}1
}
\tcbset{
  colback=black!1,boxrule=0.5pt,arc=2pt,
  left=4pt,right=4pt,top=4pt,bottom=4pt,
  fonttitle=\scriptsize\bfseries,coltitle=white,
  before skip=0pt,after skip=0pt
}
\newsavebox{\enlargementMoudafi}
\newsavebox{\enlargementBurachik}
\newsavebox{\enlargementFoundation}
\begin{lrbox}{\enlargementMoudafi}
\begin{minipage}{0.445\linewidth}
\begin{tcolorbox}[colframe=enlargeLocal,
  title={Moudafi--Oliny (2003)}]
\begin{lstlisting}[style=enlargementLean]
def enlargement
  (T : SetValuedOperator H H)
  (ε : ℝ) : SetValuedOperator H H :=
  fun x ↦ {v | ∀ y, ∀ u ∈ T y,
    -ε ≤ ⟪u - v, y - x⟫_ℝ}
\end{lstlisting}
\end{tcolorbox}
\end{minipage}
\end{lrbox}
\begin{lrbox}{\enlargementBurachik}
\begin{minipage}{0.445\linewidth}
\begin{tcolorbox}[colframe=enlargeLocal,
  title={Burachik et al. (2001)}]
\begin{lstlisting}[style=enlargementLean]
def epsilonEnlargement
  (T : SetValuedOperator H H)
  (ε : ℝ) (_hε : 0 ≤ ε) :
  SetValuedOperator H H :=
  fun x ↦ {v | ∀ z w, w ∈ T z →
    inner ℝ (w - v) (z - x) ≥ -ε}
\end{lstlisting}
\end{tcolorbox}
\end{minipage}
\end{lrbox}
\begin{lrbox}{\enlargementFoundation}
\begin{minipage}{0.445\linewidth}
\begin{tcolorbox}[colframe=enlargeShared,
  title={Foundation: proposed shared definition}]
\begin{lstlisting}[style=enlargementLean]
def enlargement
  (T : SetValuedOperator H H)
  (ε : ℝ) : SetValuedOperator H H :=
  fun x ↦ {v | ∀ y, ∀ u ∈ T y,
    -ε ≤ ⟪u - v, y - x⟫_ℝ}
\end{lstlisting}
\end{tcolorbox}
\end{minipage}
\end{lrbox}
\begin{tikzpicture}[inner sep=0pt,outer sep=0pt]
  \node (moudafi) {\usebox{\enlargementMoudafi}};
  \node[below=6pt of moudafi] (burachik)
    {\usebox{\enlargementBurachik}};
  \coordinate (mid) at ($(moudafi.east)!0.5!(burachik.east)$);
  \node[anchor=west] (foundation) at ($(mid)+(36pt,0)$)
    {\usebox{\enlargementFoundation}};
  \coordinate (join) at ($(mid)+(18pt,0)$);
  \draw[enlargeLocal,line width=0.8pt,rounded corners=3pt]
    (moudafi.east) -- ++(18pt,0) -- (join);
  \draw[enlargeLocal,line width=0.8pt,rounded corners=3pt]
    (burachik.east) -- ++(18pt,0) -- (join);
  \draw[-{Latex[length=4pt,width=3pt]},enlargeShared,line width=0.9pt]
    (join) -- (foundation.west);
\end{tikzpicture}
\endgroup
\caption{Proposed consolidation of $\varepsilon$-enlargement from
\citet{MoudafiOlny2003converge} and \citet{burachik2001robustness} into \texttt{Foundation}.
$\varepsilon\geq0$ is retained where
needed in lemmas. 
}
\label{fig:enlargement-integration}
\end{figure}

Approved candidates are consolidated into a common
definition or theorem interface, reusing an existing \texttt{Foundation} API where
appropriate. The affected paper projects are then updated and rebuilt against
the shared interface. Paper-specific algorithms and principal convergence
results remain in their respective projects.
Figure~\ref{fig:enlargement-integration} illustrates a candidate consolidation:
two paper-local definitions of $\varepsilon$-enlargement can share a single
\texttt{Foundation} interface.

Paper-level completion is checked independently of batch library integration.
After all DAG layers have reached a terminal state, \textsc{Paper2Lean} builds the complete collection of paper modules together with their import closures. A paper is marked complete only when every item is frozen, every module passes the final Lake build, and no proof placeholders remain.

\section{Evaluation}
\label{sec:ablation}

In this section, we evaluate \textsc{Paper2Lean} on 20 research papers in
mathematical optimization, focusing on inertial acceleration methods and
dynamical-system approaches. We first describe the experimental setup and
supporting libraries, and then report the scale and compilation status of
the generated projects. We next evaluate the semantic fidelity of the
formalized statements through expert review. Finally, a controlled
comparison on three representative papers examines the performance and
token consumption of two base models in statement generation and proof
construction.

\subsection{Setup}

For document extraction, statement formalization, and proof construction, we use coding agents powered by Codex CLI, with \texttt{GPT-5.6-Terra} as the base model and the reasoning effort set to \texttt{high}. All workers run as subagents.
Throughout \textsc{Paper2Lean}, all the agents interact with Lean through Lean LSP MCP~\footnote{\url{https://github.com/oOo0oOo/lean-lsp-mcp}}, which provides compiler feedback, goal-state inspection, and semantic search capabilities. All formalizations use the Lean toolchain \texttt{leanprover/lean4:v4.30.0-rc2}.

Research papers in optimization often rely on substantial background theory that is either unavailable in \texttt{mathlib} or represented through APIs that are not directly aligned with the notation and abstractions used in the source papers. To bridge this gap, we first use M2F~\citep{wang2026m2f} to formalize part of the book~\citep{bauschke2011convex} to   construct a supporting \texttt{Foundation} library containing the core definitions and results required by the target papers. We also integrate some important theorems from Optlib~\citep{li2025formalization, li2026formalization} and~\citep{bai2026formalization}. \textsc{Paper2Lean} then formalizes each paper against the combination of this \texttt{Foundation} and \texttt{mathlib}. An overview of \texttt{Foundation} is provided in Appendix~\ref{appendix:foundation}.

Table~\ref{table:1} summarizes the resulting 20 papers including source-paper pages, number of generated Lean files, lines of Lean code, token usage, and source fidelity review. In total, all projects compile successfully and contain no \texttt{sorry} placeholders.

\begin{table}[h]
\centering
\begin{tabular}{r l r r r r r}
\toprule
Number & Papers & Pages & Files & Loc & Token usage & Faithful \\ 
\midrule
1 & \citet{abbas2015dynamical} & 30 & 38 & 10819& 275,269,684 & 37/38(97\%)\\
2 & \citet{alvarez2002second} &33 & 33& 11752& 131,267,936& 27/33(82\%)\\
3 & \citet{alvarez2001inertial} & 9& 8& 1428& 11,135,901& 8/8(100\%)\\
4 & \citet{Alvarez2000minimizing} & 18& 22 & 7312& 169,776,513& 23/23(100\%)\\
5 & \citet{Alvarez2003weakconverge} &10 &14 & 4315&84,065,469 & 14/14(100\%)\\
6 &  \citet{attouch2014dynamical} & 25& 42& 6016& 176,489,473& 42/42(100\%)\\
7 & \citet{adly2023accelerated} & 43 &25 & 4892&80,166,173 & 20/25(80\%)\\
8 & \citet{ait2024second}& 26& 41& 7301&101,882,392 & 41/41(100\%)\\
9 &  \citet{attouch2024fast}& 21& 20& 4444& 84,315,239& 20/21(95\%)\\
10 &  \citet{adly2020finite}& 29& 31& 4965& 89,058,313& 28/32(88\%)\\
11 & \citet{AttouchRedont2007Newclass} &21 &25 & 8164& 241,717,003& 25/25(100\%)\\
12 & \citet{BajiCabot2006Inertial} &23 &32 & 15639& 189,536,244 & 32/32(100\%)\\
13 & \citet{FISTA2009} &20 & 20& 2210& 68,458,794& 20/20(100\%)\\
14 &  \citet{Bot2016inertialTseng}&17 & 24 & 6901 &192,617,442 & 24/24(100\%)\\
15 &  \citet{burachik2001robustness}& 20& 19& 1684& 40,900,232& 19/19(100\%)\\
16 & \citet{dong2017strong} & 22& 14& 1943& 39,369,773& 14/14(100\%)\\
17 & \citet{Guler1992proximal} & 16& 35& 2671& 95,103,480& 34/35(97\%)\\
18 &  \citet{lorenz2015inertial}& 15& 24& 4332& 58,099,932& 18/24(75\%)\\
19 & \citet{MoudafiOlny2003converge} & 8& 11& 2015& 47,541,502 & 11/11(100\%)\\
20 &  \citet{nesterov1983method}& 5& 16& 2107& 47,901,203& 16/16(100\%)\\
\bottomrule
\end{tabular}
\caption{Overview of 20 formalized papers. All generated projects compile
without \texttt{sorry}. The Faithful column reports source fidelity of the direct output of \textsc{Paper2Lean}. 
Token usage includes statement formalization and
proof construction, but excludes proof golfing.}
\label{table:1}
\end{table}

\subsection{Semantic evaluation}

Since the statement may be modified in proof stage, 
we therefore evaluate the statement by \emph{source fidelity}. 
A
statement in the output of \textsc{Paper2Lean} is called \emph{faithful} when it preserves the mathematical meaning and scope
of the source item \emph{in the context of the full paper}. The review is done after every item has been proved and every project builds without
\texttt{sorry}. This checkpoint
mainly includes statement revisions authorized during proof construction. Restoring an assumption supported by the whole source 
context is acceptable; introducing an unsupported, non-redundant assumption
that restricts the original claim is a semantic discrepancy, even if it
enables a successful proof. Appendix~\ref{appendix:statement-revisions} presents two examples: an acceptable revision that restores missing assumptions, and an unacceptable revision that assumes key facts in the proof instead of deriving them.



The human review is conducted by four experts in mathematical
optimization with experience in large-scale Lean formalization. Each expert
initially reviews five papers, then cross-audits the five papers assigned
to another reviewer. Each of the 20 papers is therefore examined by two
reviewers, and each expert reviews ten papers in total. Reviewers compare
the Lean statements with the source papers for source fidelity, and the evaluation results can be found in Table~\ref{table:1}. 

Across the 20 papers, 473 of 494 
items are faithful, giving an item-weighted source-fidelity
rate of 95.7\%. Thirteen papers have no detected discrepancy, while the
remaining seven contain 21 unfaithful statements in total. This distribution shows that \textsc{Paper2Lean} preserves the source meaning
for the large majority of statements. It also exposes a limitation:
the 21 discrepancies are found in Lean-valid developments with no proof
placeholders. In other words, successful proof construction establishes
correctness relative to the generated statement, but does not by itself
establish source fidelity.


\subsection{Model comparison}

To estimate the sensitivity of \textsc{Paper2Lean} to the base model, we conduct a controlled study on three selected papers. 
All
configurations use the same Lean, \texttt{mathlib} and  \texttt{Foundation} version,
tool access, and iteration limits. We report statement generation and proof construction separately. 

\subsubsection{Statement generation}


Unlike \emph{source fidelity}, we evaluate the quality of generated statement by \emph{input alignment}. To evaluate this index, the statements generated in formalization stage (before proof construction) are assessed against the supplied JSON items. Theorem proof bodies may still
contain \texttt{sorry}. A statement is called \emph{aligned} when it faithfully expresses the natural language in the JSON item. Two reviewers examine all three test papers.

Using the same JSON input for each paper, we compare two base models—\texttt{GPT-5.6-Terra (high)} and \texttt{DeepSeek-V4.1-Flash (high)},  both run through the Codex CLI coding agent. We record the number of statements judged acceptable by the human review in Section~\ref{sec:ablation}, together with the token usage and the number of formalize-check rounds. Table~\ref{table:stmt} reports the results.

\begin{table}[h]
\centering
\small
\begin{tabular}{ccccccc}
\toprule
\multirow{2}{*}{Papers}
    & \multicolumn{3}{c}{\texttt{GPT-5.6-Terra}}
    & \multicolumn{3}{c}{\texttt{DeepSeek-V4.1-Flash}} \\
\cmidrule(lr){2-4}
\cmidrule(lr){5-7}
    & Aligned & Token & Round
    & Aligned & Token & Round \\
\midrule
\citet{MoudafiOlny2003converge}
    & 11/11 & 13,391,172 & 1.55
    & 11/11 & 9,236,271 & 1.09 \\
\citet{abbas2015dynamical}
    & 37/38 & 60,160,380 & 1.84
    & 37/38 & 56,656,127 & 1.03 \\
\citet{BajiCabot2006Inertial}
    & 32/32 & 38,377,758 & 1.34
    & 32/32 & 26,119,725 & 1.04 \\
\bottomrule
\end{tabular}
\caption{Semantically aligned items, token usage (only statement formalization) and average formalize-check rounds in statement formalization stage under the same 15-round budget. The Aligned columns report input alignment before proof construction.}
\label{table:stmt}
\end{table}

Both models perform well in statement generation, with token consumption
on a comparable scale. \texttt{DeepSeek-V4.1-Flash} uses fewer tokens on each of the
three papers, with a modest aggregate reduction of approximately $17.8\%$
relative to \texttt{GPT-5.6-Terra}. These results suggest that both models can support
the statement-generation stage, with \texttt{DeepSeek-V4.1-Flash} offering somewhat
lower token usage in this comparison.

\subsubsection{Proof construction}

Using the same fixed statements, we compare proof construction in terms of
proved items, token consumption, and plan-prove rounds, as shown in
Table~\ref{table:prove}.

\begin{table}[h]
\centering
\small
\begin{tabular}{ccccccc}
\toprule
\multirow{2}{*}{Papers}
    & \multicolumn{3}{c}{\texttt{GPT-5.6-Terra}}
    & \multicolumn{3}{c}{\texttt{DeepSeek-V4.1-Flash}} \\
\cmidrule(lr){2-4}
\cmidrule(lr){5-7}
    & Proved & Token & Round
    & Proved & Token & Round \\
\midrule
\citet{MoudafiOlny2003converge}
    & 11/11 & 34,150,330 & 1.25
    & 11/11 & 44,088,265 & 1 \\
\citet{abbas2015dynamical}
    & 38/38 & 215,109,304 & 1.6
    & 38/38 & 926,283,566 & 1.47 \\
\citet{BajiCabot2006Inertial}
    & 32/32 & 151,158,486 & 1.94
    & 32/32 & 394,107,015 & 1.22 \\
\bottomrule
\end{tabular}
\caption{Proved items, token usage, and average plan-prove rounds in proof construction under the same 10-round budget. Token usage includes only proof construction, excluding proof golfing.}
\label{table:prove}
\end{table}

Both models complete all 81 proof-construction tasks across the three papers
within the same iteration budget, including demanding items whose Lean
developments exceed 1,000 lines. However, \texttt{DeepSeek-V4.1-Flash} consumes
substantially more tokens: its total usage is approximately $3.41$ times
that of \texttt{GPT-5.6-Terra}, with per-paper ratios ranging from $1.29$ to $4.31$.
Thus, both models can handle these proof tasks, but \texttt{GPT-5.6-Terra} achieves
the same completion rate with substantially lower token usage under the
evaluated workflow.

\subsection{OptPaper-20: dataset after human refinement}

For all 20 papers formalized by \textsc{Paper2Lean}, we manually refine the generated code, primarily to improve the quality of the statements and remove redundant definitions or assumptions. These changes make the code more readable and reusable. Examples of human refinement can be found in the next section. We also translate the refined formal statements into natural language. The complete \textsc{OptPaper-20} dataset contains the Lean files for all 20 papers, their corresponding natural language descriptions, and structured JSON files.

\section{Case studies and discussions}

In this section, we examine two cases of \textsc{Paper2Lean} output. The first is semantically faithful in the context of the full paper but can be refined to make downstream theorem statements more concise and readable. The second passes the agent's semantic
check and Lean proof verification yet is judged unfaithful because it
omits the source claim and its contextual dependencies.

\subsection{Case 1: faithful, but requiring refinement.}

We judge the original trajectory definition for the Hessian-driven inertial
equation (DIN) in
Figure~\ref{fig:statement-comparison} \emph{semantically faithful} in the
context of \citet{alvarez2002second}: downstream assumptions identify the
supplied gradient and Hessian maps with those of the potential, while
\texttt{HasDerivAt} conditions identify velocity and acceleration. However,
these auxiliary parameters make the interface cumbersome to reuse. The
refinement expresses derivatives directly through \texttt{mathlib}'s
\texttt{gradient}, \texttt{fderiv}, and \texttt{deriv}, and restricts the
equation to $t\geq 0$, matching its use in subsequent theorems.
Appendix~\ref{appendix:din-downstream} presents the corresponding downstream
assumptions and weak convergence theorem before and after refinement.
Removing auxiliary map parameters and repeated compatibility hypotheses
makes the downstream theorem statements less verbose and easier to read.

\begin{figure}[htbp]
\centering
\begingroup
\definecolor{stmtBefore}{HTML}{536A82}
\definecolor{stmtAfter}{HTML}{287B70}
\lstdefinestyle{leanStatementComparison}{
  basicstyle=\ttfamily\scriptsize,
  keywordstyle=\color{blue!55!black}\bfseries,
  commentstyle=\color{black!45},
  morekeywords={def,Prop},
  literate={ℝ}{{$\mathbb{R}$}}1
           {Φ}{{$\Phi$}}1
           {α}{{$\alpha$}}1
           {β}{{$\beta$}}1
           {→}{{$\rightarrow$}}1
           {∀}{{$\forall$}}1
           {∧}{{$\land$}}1
           {•}{{$\bullet$}}1
           {≤}{{$\leq$}}1,
  morecomment=[l]{--},
  columns=fullflexible,
  keepspaces=true,
  showstringspaces=false,
  breaklines=true,
  aboveskip=0pt,
  belowskip=0pt
}
\tcbset{
  colback=black!1,
  boxrule=0.6pt,
  arc=2pt,
  valign=top,
  equal height group=statementComparison,
  left=7pt,right=7pt,top=3pt,bottom=3pt,
  fonttitle=\small\bfseries,
  coltitle=white,
  before skip=0pt,after skip=0pt
}
\begin{minipage}[t]{0.485\linewidth}
\vspace{0pt}
\begin{tcolorbox}[colframe=stmtBefore,
  title={(a) Before: faithful statement}]
\begin{lstlisting}[style=leanStatementComparison]
def IsDINTrajectory (H : Type*)
  [NormedAddCommGroup H]
  [InnerProductSpace ℝ H]
  (gradΦ : H → H) 
  (hessΦ : H → H →L[ℝ] H)
  (α β : ℝ) (x v a : ℝ → H) : Prop :=
  0 < α ∧ 0 < β ∧
  (∀ t, HasDerivAt x (v t) t) ∧
  (∀ t, HasDerivAt v (a t) t) ∧
  ∀ t, a t + α • v t +
  β • hessΦ (x t) (v t) + 
  gradΦ (x t) = 0
\end{lstlisting}
\end{tcolorbox}
\end{minipage}\hfill
\begin{minipage}[t]{0.485\linewidth}
\vspace{0pt}
\begin{tcolorbox}[colframe=stmtAfter,
  title={(b) After: refined statement}]
\begin{lstlisting}[style=leanStatementComparison]
def IsDINTrajectory (H : Type*)
  [NormedAddCommGroup H]
  [InnerProductSpace ℝ H]
  [CompleteSpace H] (Φ : H → ℝ)
  (α β : ℝ) (x : ℝ → H) : Prop :=
  0 < α ∧ 0 < β ∧ ∀ t, 0 ≤ t →
  deriv (deriv x) t + α • deriv x t +
  β • (fderiv ℝ (gradient Φ) (x t))
  (deriv x t) + gradient Φ (x t) = 0
\end{lstlisting}
\end{tcolorbox}
\end{minipage}
\endgroup
\caption{DIN trajectory in \citet{alvarez2002second}:
(a) original statement; (b) refined statement.}
\label{fig:statement-comparison}
\end{figure}

This refinement improves usability but requires care with assumptions.
The predicates are not unconditionally equivalent because their time domains
differ, and downstream results must retain the required regularity of $\Phi$
and twice differentiability of the trajectory: \texttt{deriv} alone does not
assert differentiability. The case therefore motivates refining a semantically faithful statement for reuse while preserving its mathematical obligations.

\subsection{Case 2: unfaithful despite passing agent checks and proof.}
The formalization of Remark~2.19 in \citet{abbas2015dynamical}
(Figure~\ref{fig:kadec-klee-failure}) passed the formalization agent's semantic check and
was successfully proved in Lean, but we judge it \emph{not semantically
faithful}. The source claims strong convergence of the system's trajectory
under the Kadec--Klee property. Under the standing assumptions, including
$S\neq\emptyset$, Theorem~2.8 and Proposition~2.16 provide weak convergence
to $x_\infty\in S$ and convergence of $\Phi(x(t))$ to $\Phi(x_\infty)$;
applying the property along sequences $t_k\to+\infty$ yields
$\|x(t)-x_\infty\|\to 0$.

\begin{figure}[htbp]
\centering
\begingroup
\lstdefinestyle{leanKadecKleeFailure}{
  basicstyle=\ttfamily\scriptsize,
  keywordstyle=\color{blue!55!black}\bfseries,
  morekeywords={theorem,Type,fun},
  literate={ℝ}{{$\mathbb{R}$}}1
           {ℕ}{{$\mathbb{N}$}}1
           {→}{{$\rightarrow$}}1
           {↦}{{$\mapsto$}}1
           {∀}{{$\forall$}}1
           {∈}{{$\in$}}1
           {⊥}{{$\bot$}}1
           {‖}{{$\Vert$}}1,
  columns=fullflexible,
  keepspaces=true,
  showstringspaces=false,
  breaklines=true,
  aboveskip=2pt,
  belowskip=2pt
}
\begin{tcolorbox}[
  colback=black!1,colframe=black!50,
  boxrule=0.6pt,arc=2pt,
  left=7pt,right=7pt,top=7pt,bottom=7pt,
  fonttitle=\small\bfseries,
  title={(a) Natural-language source: Remark 2.19}]
\small
``By a similar argument, strong convergence of $x(\cdot)$ holds under the
Kadec-Klee property: whenever $(x_k)$ converges weakly to $x$, and
$(\Phi(x_k))$ converges to $\Phi(x)$, then $(x_k)$ converges strongly to $x$.''
\end{tcolorbox}
\begin{tcolorbox}[
  colback=black!1,colframe=black!50,
  boxrule=0.6pt,arc=2pt,
  left=7pt,right=7pt,top=7pt,bottom=7pt,
  fonttitle=\small\bfseries,
  title={(b) Generated Lean statement and proof}]
\begin{lstlisting}[style=leanKadecKleeFailure]
theorem kadec_klee_property {H : Type*}
    [NormedAddCommGroup H] [NormedSpace ℝ H] [InnerProductSpace ℝ H]
    (Phi : H → Set.Ioi (⊥ : EReal))
    (_hPhi : Phi ∈ ERealFunction.gammaZero H) :
    (hKadecKlee : ∀ (xSeq : ℕ → H) (x : H),
      Tendsto (fun k ↦ toWeakSpace ℝ H (xSeq k)) atTop (nhds (toWeakSpace ℝ H x)) →
      Tendsto (fun k ↦ (Phi (xSeq k) : EReal)) atTop (nhds (Phi x : EReal)) →
      Tendsto (fun k ↦ ‖xSeq k - x‖) atTop (nhds 0)) → ∀ (xSeq : ℕ → H) (x : H),
      Tendsto (fun k ↦ toWeakSpace ℝ H (xSeq k)) atTop (nhds (toWeakSpace ℝ H x)) →
      Tendsto (fun k ↦ (Phi (xSeq k) : EReal)) atTop (nhds (Phi x : EReal)) →
      Tendsto (fun k ↦ ‖xSeq k - x‖) atTop (nhds 0) := id
\end{lstlisting}
\end{tcolorbox}
\endgroup
\caption{Remark~2.19 of \citet{abbas2015dynamical} and its generated
formalization, which passed automatic semantic checking and Lean proof
verification.}
\label{fig:kadec-klee-failure}
\end{figure}

The hypothesis \texttt{hKadecKlee} and the conclusion
express the same property, so \texttt{id} proves the statement. The claim
about strong convergence of the source paper's trajectory is absent.
The assumption itself is legitimate, but omitting its consequence for the trajectory loses the source claim. This illustrates how an agent can miss
the context behind ``by a similar argument'' even when semantic checking
and proof construction both succeed. Review must therefore check the
intended conclusion against the surrounding paper and its dependencies.

\section{Conclusion}






We introduced \textsc{Paper2Lean}, an agentic framework for converting
mathematical research papers into integrated Lean projects through
dependency-aware extraction, formalization, semantic checking, theorem proving,
and library integration. Its central representation is a dependency DAG that
supports staged processing, compilation-gated validation, and the reuse of
verified interfaces across a complete paper. Evaluated on 20 optimization
papers with \texttt{mathlib} and a shared supporting library, the framework
generated projects that compile without proof placeholders. The resulting
\textsc{OptPaper-20} corpus contains 494 Lean files and 110,910 lines of code,
together with the natural language and structured JSON files.

Our audit also shows that Lean-checked compilation and faithful
formalization measure different properties. Before manual correction, 473 of
494 reviewed declarations were judged faithful; all identified
discrepancies were repaired in the released corpus. This finding
supports paper-level formalization as a practical direction while showing that
semantic review remains necessary for detecting omitted assumptions, altered
claims, and lost contextual dependencies. Future work will focus on automated
semantic checking~\citep{liu2026assess,liu2025rethinking}, better recovery of implicit assumptions, lower token costs,
broader mathematical domains, and more systematic comparisons with alternative
agent architectures and formalization pipelines.

\bibliography{references}

\begin{thebibliography}{47}
\providecommand{\natexlab}[1]{#1}
\providecommand{\url}[1]{\texttt{#1}}
\expandafter\ifx\csname urlstyle\endcsname\relax
  \providecommand{\doi}[1]{doi: #1}\else
  \providecommand{\doi}{doi: \begingroup \urlstyle{rm}\Url}\fi

\bibitem[Abbas \& Attouch(2015)Abbas and Attouch]{abbas2015dynamical}
Boushra Abbas and Hedy Attouch.
\newblock Dynamical systems and forward--backward algorithms associated with
  the sum of a convex subdifferential and a monotone cocoercive operator.
\newblock \emph{Optimization}, 64\penalty0 (10):\penalty0 2223--2252, 2015.

\bibitem[Adly \& Attouch(2020)Adly and Attouch]{adly2020finite}
Samir Adly and Hedy Attouch.
\newblock Finite convergence of proximal-gradient inertial algorithms combining
  dry friction with hessian-driven damping.
\newblock \emph{SIAM Journal on Optimization}, 30\penalty0 (3):\penalty0
  2134--2162, 2020.

\bibitem[Adly \& Attouch(2023)Adly and Attouch]{adly2023accelerated}
Samir Adly and Hedy Attouch.
\newblock Accelerated dynamics with dry friction via time scaling and averaging
  of doubly nonlinear evolution equations.
\newblock \emph{Nonlinear Analysis: Hybrid Systems}, 50:\penalty0 101402, 2023.

\bibitem[Ait~Mansour et~al.(2024)Ait~Mansour, Attouch, and
  Mazgouri]{ait2024second}
Mohamed Ait~Mansour, Hedy Attouch, and Zakaria Mazgouri.
\newblock Second order dynamical systems governed by maximally monotone
  operators: Quantitative stability.
\newblock \emph{Evol. Equ. Control Theory}, 13:\penalty0 510--535, 2024.

\bibitem[Alvarez(2000)]{Alvarez2000minimizing}
Felipe Alvarez.
\newblock On the minimizing property of a second order dissipative system in
  hilbert spaces.
\newblock \emph{SIAM Journal on Control and Optimization}, 38\penalty0
  (4):\penalty0 1102--1119, 2000.

\bibitem[Alvarez(2003)]{Alvarez2003weakconverge}
Felipe Alvarez.
\newblock Weak convergence of a relaxed and inertial hybrid projection-proximal
  point algorithm for maximal monotone operators in hilbert space.
\newblock \emph{SIAM J. on Optimization}, 14\penalty0 (3):\penalty0 773–782,
  March 2003.

\bibitem[Alvarez \& Attouch(2001)Alvarez and Attouch]{alvarez2001inertial}
Felipe Alvarez and Hedy Attouch.
\newblock An inertial proximal method for maximal monotone operators via
  discretization of a nonlinear oscillator with damping.
\newblock \emph{Set-Valued Analysis}, 9\penalty0 (1):\penalty0 3--11, 2001.

\bibitem[Alvarez et~al.(2002)Alvarez, Attouch, Bolte, and
  Redont]{alvarez2002second}
Felipe Alvarez, Hedy Attouch, J{\'e}r{\^o}me Bolte, and Patrick Redont.
\newblock A second-order gradient-like dissipative dynamical system with
  hessian-driven damping.: Application to optimization and mechanics.
\newblock \emph{Journal de math{\'e}matiques pures et appliqu{\'e}es},
  81\penalty0 (8):\penalty0 747--779, 2002.

\bibitem[{Anthropic}(2026)]{anthropic2026flt}
{Anthropic}.
\newblock Fermat's last theorem in {Lean 4}.
\newblock \url{https://github.com/anthropics/fermats-last-theorem}, September
  2026.
\newblock Complete machine-checked formalization; accessed September 22, 2026.

\bibitem[Attouch et~al.(2007)Attouch, Redont, and
  Soubeyran]{AttouchRedont2007Newclass}
H.~Attouch, P.~Redont, and A.~Soubeyran.
\newblock A new class of alternating proximal minimization algorithms with
  costs-to-move.
\newblock \emph{SIAM Journal on Optimization}, 18\penalty0 (3):\penalty0
  1061--1081, 2007.

\bibitem[Attouch et~al.(2014)Attouch, Peypouquet, and
  Redont]{attouch2014dynamical}
H{\'e}dy Attouch, Juan Peypouquet, and Patrick Redont.
\newblock A dynamical approach to an inertial forward-backward algorithm for
  convex minimization.
\newblock \emph{SIAM Journal on Optimization}, 24\penalty0 (1):\penalty0
  232--256, 2014.

\bibitem[Attouch et~al.(2024)Attouch, Chbani, and Riahi]{attouch2024fast}
Hedy Attouch, Zaki Chbani, and Hassan Riahi.
\newblock Fast convex optimization via a third-order in time evolution
  equation: Toges-v an improved version of toges.
\newblock \emph{Optimization}, 73\penalty0 (3):\penalty0 575--595, 2024.

\bibitem[Azerbayev et~al.(2023)Azerbayev, Piotrowski, Schoelkopf, Ayers, Radev,
  and Avigad]{azerbayev2023proofnet}
Zhangir Azerbayev, Bartosz Piotrowski, Hailey Schoelkopf, Edward~W Ayers,
  Dragomir Radev, and Jeremy Avigad.
\newblock Proofnet: Autoformalizing and formally proving undergraduate-level
  mathematics.
\newblock \emph{arXiv preprint arXiv:2302.12433}, 2023.

\bibitem[Bai et~al.(2026)Bai, Li, Yu, and Liang]{bai2026formalization}
Yifan Bai, Yantao Li, Jian Yu, and Jingwei Liang.
\newblock Formalization of two fixed-point algorithms in hilbert spaces.
\newblock \emph{arXiv preprint arXiv:2602.17064}, 2026.

\bibitem[Baji \& Cabot(2006)Baji and Cabot]{BajiCabot2006Inertial}
Bruno Baji and Alexandre Cabot.
\newblock An inertial proximal algorithm with dry friction: Finite convergence
  results.
\newblock \emph{Set-Valued Analysis}, 14:\penalty0 1--23, 01 2006.

\bibitem[Barras et~al.(1997)Barras, Boutin, Cornes, Courant, Filliatre,
  Gimenez, Herbelin, Huet, Munoz, Murthy, et~al.]{barras1997coq}
Bruno Barras, Samuel Boutin, Cristina Cornes, Judica{\"e}l Courant,
  Jean-Christophe Filliatre, Eduardo Gimenez, Hugo Herbelin, Gerard Huet, Cesar
  Munoz, Chetan Murthy, et~al.
\newblock \emph{The Coq proof assistant reference manual: Version 6.1}.
\newblock PhD thesis, Inria, 1997.

\bibitem[Bauschke \& Combettes(2011)Bauschke and Combettes]{bauschke2011convex}
H.~Bauschke and P.~L. Combettes.
\newblock \emph{Convex Analysis and Monotone Operator Theory in Hilbert
  Spaces}.
\newblock Springer, 2011.

\bibitem[Beck \& Teboulle(2009)Beck and Teboulle]{FISTA2009}
Amir Beck and Marc Teboulle.
\newblock A fast iterative shrinkage-thresholding algorithm for linear inverse
  problems.
\newblock \emph{SIAM J. Img. Sci.}, 2\penalty0 (1):\penalty0 183–202, March
  2009.

\bibitem[Bo{\'z} \& Csetnek(2016)Bo{\'z} and Csetnek]{Bot2016inertialTseng}
Radu~Ioan Bo{\'z} and Ern{\"o}~Robert Csetnek.
\newblock An inertial tseng's type proximal algorithm for nonsmooth and
  nonconvex optimization problems.
\newblock \emph{J. Optim. Theory Appl.}, 171\penalty0 (2):\penalty0 600–616,
  November 2016.

\bibitem[Burachik et~al.(2001)Burachik, Scheimberg, and
  Svaiter]{burachik2001robustness}
Regina~S Burachik, Susana Scheimberg, and Benar~F Svaiter.
\newblock Robustness of the hybrid extragradient proximal-point algorithm.
\newblock \emph{Journal of Optimization Theory and Applications}, 111\penalty0
  (1):\penalty0 117--136, 2001.

\bibitem[Cabral et~al.(2026)Cabral, Manh, Yu, Tai, Feng, and
  Xin]{cabral2026proofflow}
Rafael Cabral, Tuan Manh, Xuejun Yu, Wai~Ming Tai, Zijin Feng, and Shen Xin.
\newblock Proofflow: A dependency graph approach to faithful proof
  autoformalization.
\newblock In \emph{International Conference on Learning Representations},
  volume 2026, pp.\  58459--58483, 2026.

\bibitem[De~Moura et~al.(2015)De~Moura, Kong, Avigad, Van~Doorn, and von
  Raumer]{de2015lean}
Leonardo De~Moura, Soonho Kong, Jeremy Avigad, Floris Van~Doorn, and Jakob von
  Raumer.
\newblock The lean theorem prover (system description).
\newblock In \emph{International Conference on Automated Deduction}, pp.\
  378--388. Springer, 2015.

\bibitem[DeepSeek-AI(2026)]{deepseek-harness2026}
DeepSeek-AI.
\newblock Deepseek harness: Everything is a plugin.
\newblock \url{https://github.com/deepseek-ai/deepseek-harness}, 2026.

\bibitem[Dong et~al.(2017)Dong, Jiang, Cholamjiak, and Shehu]{dong2017strong}
Qiao-Li Dong, Dan Jiang, Prasit Cholamjiak, and Yekini Shehu.
\newblock A strong convergence result involving an inertial forward–backward
  algorithm for monotone inclusions.
\newblock \emph{Journal of Fixed Point Theory and Applications}, 19, 12 2017.

\bibitem[Guler(1992)]{Guler1992proximal}
Osman Guler.
\newblock New proximal point algorithms for convex minimization.
\newblock \emph{SIAM Journal on Optimization}, 2:\penalty0 649--664, 01 1992.

\bibitem[Hanyu et~al.(2026)Hanyu, Xie, Yutong, Gao, and Dong]{hanyu2026aria}
Wang Hanyu, Ruohan Xie, Wang Yutong, Guoxiong Gao, and Bin Dong.
\newblock Aria: An agent for retrieval and iterative auto-formalization via
  dependency graph.
\newblock In \emph{International Conference on Learning Representations},
  volume 2026, pp.\  11690--11722, 2026.

\bibitem[Ju et~al.(2026{\natexlab{a}})Ju, Chen, Wu, Dai, and
  Dong]{ju2026matlas}
Haocheng Ju, Leheng Chen, Peihao Wu, Bryan Dai, and Bin Dong.
\newblock Matlas: A semantic search engine for mathematics.
\newblock \emph{arXiv preprint arXiv:2604.17484}, 2026{\natexlab{a}}.

\bibitem[Ju et~al.(2026{\natexlab{b}})Ju, Gao, Jiang, Wu, Sun, Liu, Chen, Wang,
  Wang, Wang, et~al.]{ju2026automated}
Haocheng Ju, Guoxiong Gao, Jiedong Jiang, Bin Wu, Zeming Sun, Shurui Liu,
  Leheng Chen, Yutong Wang, Yuefeng Wang, Zichen Wang, et~al.
\newblock Automated conjecture resolution with formal verification.
\newblock \emph{arXiv preprint arXiv:2604.03789}, 2026{\natexlab{b}}.

\bibitem[Li et~al.(2025)Li, Wang, He, Wu, Xu, and Wen]{li2025formalization}
Chenyi Li, Ziyu Wang, Wanyi He, Yuxuan Wu, Shengyang Xu, and Zaiwen Wen.
\newblock Formalization of convergence rates of four first-order algorithms for
  convex optimization: C. li et al.
\newblock \emph{Journal of Automated Reasoning}, 69\penalty0 (4):\penalty0 28,
  2025.

\bibitem[Li et~al.(2026)Li, Wang, Bai, Duan, Gao, Hao, and
  Wen]{li2026formalization}
Chenyi Li, Zichen Wang, Yifan Bai, Yunxi Duan, Yuqing Gao, Pengfei Hao, and
  Zaiwen Wen.
\newblock Formalization of algorithms for optimization with block structures.
\newblock \emph{Science China Mathematics}, pp.\  1--18, 2026.

\bibitem[Liu et~al.(2026{\natexlab{a}})Liu, Zhou, Zhu, Santos, He, Liu, Wang,
  Xie, Zhao, Wang, et~al.]{liu2026numina}
Junqi Liu, Zihao Zhou, Zekai Zhu, Marco~Dos Santos, Weikun He, Jiawei Liu, Ran
  Wang, Yunzhou Xie, Junqiao Zhao, Qiufeng Wang, et~al.
\newblock Numina-lean-agent: An open and general agentic reasoning system for
  formal mathematics.
\newblock \emph{arXiv preprint arXiv:2601.14027}, 2026{\natexlab{a}}.

\bibitem[Liu et~al.(2025)Liu, Zheng, Lu, Cao, and Yan]{liu2025rethinking}
Qi~Liu, Xinhao Zheng, Xudong Lu, Qinxiang Cao, and Junchi Yan.
\newblock Rethinking and improving autoformalization: Towards a faithful metric
  and a dependency retrieval-based approach.
\newblock In \emph{International Conference on Learning Representations},
  volume 2025, pp.\  86205--86233, 2025.

\bibitem[Liu et~al.(2026{\natexlab{b}})Liu, Bao, Zhang, Liu, Chen, Liu, Jiao,
  and Luo]{liu2026atlas}
Xiaoyang Liu, Kangjie Bao, Jiashuo Zhang, Yunqi Liu, Yu~Chen, Yuntian Liu, Yang
  Jiao, and Tao Luo.
\newblock Atlas: Autoformalizing theorems through lifting, augmentation, and
  synthesis of data.
\newblock \emph{Advances in Neural Information Processing Systems},
  38:\penalty0 56290--56326, 2026{\natexlab{b}}.

\bibitem[Liu et~al.(2026{\natexlab{c}})Liu, Dong, Bai, Li, Liu, and
  Luo]{liu2026decompose}
Xiaoyang Liu, Zineng Dong, Yifan Bai, Yantao Li, Yuntian Liu, and Tao Luo.
\newblock Decompose, structure, and repair: A neuro-symbolic framework for
  autoformalization via operator trees.
\newblock \emph{arXiv preprint arXiv:2604.19000}, 2026{\natexlab{c}}.

\bibitem[Liu et~al.(2026{\natexlab{d}})Liu, Zhu, Dong, Liu, ZhaoXuan, Chen,
  Luo, et~al.]{liu2026assess}
Xiaoyang Liu, Tao Zhu, Zineng Dong, Yuntian Liu, Liu ZhaoXuan, Yu~Chen, Tao
  Luo, et~al.
\newblock Assess: A semantic and structural evaluation framework for statement
  similarity.
\newblock In \emph{International Conference on Learning Representations},
  volume 2026, pp.\  143709--143734, 2026{\natexlab{d}}.

\bibitem[Lorenz \& Pock(2015)Lorenz and Pock]{lorenz2015inertial}
Dirk~A Lorenz and Thomas Pock.
\newblock An inertial forward-backward algorithm for monotone inclusions.
\newblock \emph{Journal of Mathematical Imaging and Vision}, 51\penalty0
  (2):\penalty0 311--325, 2015.

\bibitem[mathlib {C}ommunity(2020)]{mathlib2020}
The mathlib {C}ommunity.
\newblock The {L}ean {M}athematical {L}ibrary.
\newblock In \emph{Proceedings of the 9th {ACM} {SIGPLAN} {I}nternational
  {C}onference on {C}ertified {P}rograms and {P}roofs}, {CPP} '20. {ACM},
  {J}anuary 2020.

\bibitem[Moudafi \& Oliny(2003)Moudafi and Oliny]{MoudafiOlny2003converge}
A.~Moudafi and M.~Oliny.
\newblock Convergence of a splitting inertial proximal method for monotone
  operators.
\newblock \emph{Journal of Computational and Applied Mathematics}, 155\penalty0
  (2):\penalty0 447--454, 2003.

\bibitem[Nesterov(1983)]{nesterov1983method}
Yurii Nesterov.
\newblock A method for solving the convex programming problem with convergence
  rate o (1/k2).
\newblock In \emph{Dokl akad nauk Sssr}, volume 269, pp.\  543, 1983.

\bibitem[Nipkow et~al.(2002)Nipkow, Wenzel, and Paulson]{nipkow2002isabelle}
Tobias Nipkow, Markus Wenzel, and Lawrence~C Paulson.
\newblock \emph{Isabelle/HOL: a proof assistant for higher-order logic}.
\newblock Springer, 2002.

\bibitem[Ren et~al.(2025)Ren, Shao, Song, Xin, Wang, Zhao, Zhang, Fu, Zhu,
  Yang, et~al.]{ren2025deepseek}
ZZ~Ren, Zhihong Shao, Junxiao Song, Huajian Xin, Haocheng Wang, Wanjia Zhao,
  Liyue Zhang, Zhe Fu, Qihao Zhu, Dejian Yang, et~al.
\newblock Deepseek-prover-v2: Advancing formal mathematical reasoning via
  reinforcement learning for subgoal decomposition.
\newblock \emph{arXiv preprint arXiv:2504.21801}, 2025.

\bibitem[Tsoukalas et~al.(2026)Tsoukalas, Kovsharov, Shirobokov, Surina,
  Firsching, B{\'e}rczi, Ruiz, Suggala, Wagner, Wieser, Yu, Huang, Horv{\'a}th,
  Ferrauiolo, Michalewski, Lockhart, Grosu, Hubert, Balog, Kohli, and
  Chaudhuri]{tsoukalas2026advancing}
George Tsoukalas, Anton Kovsharov, Sergey Shirobokov, Anja Surina, Moritz
  Firsching, Gergely B{\'e}rczi, Francisco J.~R. Ruiz, Arun Suggala, Adam~Zsolt
  Wagner, Eric Wieser, Lei Yu, Aja Huang, Mikl{\'o}s~Z. Horv{\'a}th, Andrew
  Ferrauiolo, Henryk Michalewski, Edward Lockhart, Codrut Grosu, Thomas Hubert,
  Matej Balog, Pushmeet Kohli, and Swarat Chaudhuri.
\newblock Advancing mathematics research with {AI}-driven formal proof search.
\newblock \emph{arXiv preprint arXiv:2605.22763}, 2026.

\bibitem[Wang et~al.(2024)Wang, Xu, Zhao, Ouyang, Wu, Zhao, Xu, Liu, Qu, Shang,
  et~al.]{wang2024mineru}
Bin Wang, Chao Xu, Xiaomeng Zhao, Linke Ouyang, Fan Wu, Zhiyuan Zhao, Rui Xu,
  Kaiwen Liu, Yuan Qu, Fukai Shang, et~al.
\newblock Mineru: An open-source solution for precise document content
  extraction.
\newblock \emph{arXiv preprint arXiv:2409.18839}, 2024.

\bibitem[Wang et~al.(2025)Wang, Unsal, Lin, Baksys, Liu, Santos, Sung, Vinyes,
  Ying, Zhu, et~al.]{wang2025kimina}
Haiming Wang, Mert Unsal, Xiaohan Lin, Mantas Baksys, Junqi Liu, Marco~Dos
  Santos, Flood Sung, Marina Vinyes, Zhenzhe Ying, Zekai Zhu, et~al.
\newblock Kimina-prover preview: Towards large formal reasoning models with
  reinforcement learning.
\newblock \emph{arXiv preprint arXiv:2504.11354}, 2025.

\bibitem[Wang et~al.(2026)Wang, Ma, Ming, Zhang, Yuan, and Wen]{wang2026m2f}
Zichen Wang, Wanli Ma, Zhenyu Ming, Gong Zhang, Kun Yuan, and Zaiwen Wen.
\newblock M2f: Automated formalization of mathematical literature at scale.
\newblock \emph{arXiv preprint arXiv:2602.17016}, 2026.

\bibitem[Yang et~al.(2023)Yang, Swope, Gu, Chalamala, Song, Yu, Godil, Prenger,
  and Anandkumar]{yang2023leandojo}
Kaiyu Yang, Aidan Swope, Alex Gu, Rahul Chalamala, Peiyang Song, Shixing Yu,
  Saad Godil, Ryan~J Prenger, and Animashree Anandkumar.
\newblock Leandojo: Theorem proving with retrieval-augmented language models.
\newblock \emph{Advances in Neural Information Processing Systems},
  36:\penalty0 21573--21612, 2023.

\bibitem[Zheng et~al.(2021)Zheng, Han, and Polu]{zheng2021minif2f}
Kunhao Zheng, Jesse~Michael Han, and Stanislas Polu.
\newblock Minif2f: a cross-system benchmark for formal olympiad-level
  mathematics.
\newblock \emph{arXiv preprint arXiv:2109.00110}, 2021.

\end{thebibliography}
\bibliographystyle{references}

\clearpage
\appendix
\section{Example Intermediate Artifact}
\label{appendix:intermediate-artifact}

\paragraph{Extraction output.}
After document extraction, \textsc{Paper2Lean} stores the mathematical content
of each source paper in a structured JSON artifact, which mainly records paper-level metadata and a list of mathematical items.
Each item has a stable identifier and type, a human-readable name, the
extracted mathematical statement, a natural language proof (optional), and explicit
\texttt{depends\_on} edges.  These edges form the initial dependency graph
used to schedule downstream statement formalization and proof construction.
Listing~\ref{lst:fista-extraction} is a faithful excerpt from the extraction
output for \cite{FISTA2009}.  It follows a representative chain
from the optimization model and proximal map, through an optimality lemma and
the FISTA algorithm, to the paper's main convergence theorem.  All displayed
values are reproduced from the generated artifact.  For readability, the long
\texttt{proof} field of the final theorem is omitted from the
selected records.

\begin{lstlisting}[
  style=jsonArtifact,
  caption={Excerpt from \texttt{paper\_expanded.json} for \protect\citet{FISTA2009}, showing extracted definitions, a lemma, an algorithm, a theorem, and their dependency edges.},
  label={lst:fista-extraction}
]
{
  "metadata": {
    "title": "A Fast Iterative Shrinkage-Thresholding Algorithm for Linear Inverse Problems",
    "authors": [
      "Amir Beck",
      "Marc Teboulle"
    ],
    "year": 2009,
    "source_type": "markdown",
    "venue": "SIAM Journal on Imaging Sciences, Vol. 2, No. 1, pp. 183-202"
  },
  "items": [
    {
      "id": "def_5",
      "type": "definition",
      "name": "General convex optimization model (Problem P)",
      "statement": "The general nonsmooth convex optimization model is $$(\\mathrm{P})\\quad \\min\\{F(x) \\equiv f(x) + g(x) : x \\in \\mathbb{R}^n\\},$$ where:\\n- $g : \\mathbb{R}^n \\to \\mathbb{R}$ is a continuous convex function which is possibly nonsmooth.\\n- $f : \\mathbb{R}^n \\to \\mathbb{R}$ is a smooth convex function of type $C^{1,1}$, i.e., continuously differentiable with Lipschitz continuous gradient $L(f)$.\\n- Problem (P) is solvable, i.e., $X^* := \\arg\\min F \\ne \\emptyset$, and for $x^* \\in X^*$ we set $F_* := F(x^*)$.",
      "proof": null,
      "depends_on": ["def_3"]
    },
    {
      "id": "def_7",
      "type": "definition",
      "name": "Proximal operator $p_L(y)$",
      "statement": "For any $L > 0$ and $y \\in \\mathbb{R}^n$, the proximal operator (or prox-map) is defined as the unique minimizer $$p_L(y) := \\arg\\min\\{Q_L(x,y) : x \\in \\mathbb{R}^n\\}.$$ Equivalently, ignoring constant terms in $y$, $$p_L(y) = \\arg\\min_x \\left\\{ g(x) + \\frac{L}{2}\\left\\|x - \\left(y - \\frac{1}{L}\\nabla f(y)\\right)\\right\\|^2 \\right\\}.$$",
      "proof": null,
      "depends_on": ["def_6"]
    },
    {
      "id": "lem_2",
      "type": "lemma",
      "name": "Lemma 2.2: Optimality characterization of the proximal operator",
      "statement": "For any $y \\in \\mathbb{R}^n$, one has $z = p_L(y)$ if and only if there exists $\\gamma(y) \\in \\partial g(z)$, the subdifferential of $g(\\cdot)$, such that $$\\nabla f(y) + L(z-y) + \\gamma(y) = 0.$$",
      "proof": "The proof is immediate from optimality conditions for the strongly convex problem (2.6). Since $Q_L(\\cdot, y)$ is strongly convex with modulus $L$, the first-order necessary and sufficient optimality condition is $0 \\in \\partial(Q_L(\\cdot, y))(z) = \\nabla f(y) + L(z-y) + \\partial g(z)$.",
      "depends_on": [
        "def_4",
        "def_7"
      ]
    },
    {
      "id": "algo_1",
      "type": "algorithm",
      "name": "FISTA with constant stepsize",
      "statement": "FISTA with constant stepsize: the sequences satisfy the proximal update and the momentum recurrences of equations (4.1)--(4.3).",
      "proof": null,
      "depends_on": ["def_5", "def_7"]
    },
    {
      "id": "thm_2",
      "type": "theorem",
      "name": "Theorem 4.4: Global convergence rate of FISTA -- $O(1/k^2)$",
      "statement": "Let $\\{x_k\\}$, $\\{y_k\\}$ be generated by FISTA with either constant or backtracking stepsize rule. Then for any $k \\ge 1$, $$F(x_k) - F(x^*) \\le \\frac{2\\alpha L(f)\\|x_0 - x^*\\|^2}{(k+1)^2} \\quad \\forall x^* \\in X^*,$$ where $\\alpha = 1$ for the constant stepsize setting and $\\alpha = \\eta$ for the backtracking stepsize setting.",
      "depends_on": [
        "def_5",
        "def_7",
        "lem_3",
        "lem_4",
        "lem_5",
        "lem_6"
      ]
    }
  ]
}
\end{lstlisting}

The excerpt illustrates how heterogeneous paper content is represented in one
uniform schema.  The definitions introduce the composite objective and its
proximal map; \texttt{lem\_2} records the first-order optimality
characterization of that map; \texttt{algo\_1} identifies the FISTA update
scheme; and \texttt{thm\_2} states the accelerated $O(1/k^2)$ objective-value
rate.  The explicit dependency lists retain the semantic prerequisites for
each item.  In particular, the theorem depends on the model and proximal-map
definitions as well as the four lemmas used in its proof, allowing downstream
agents to formalize and verify these declarations in dependency order.

\section{Shared Foundation Library}
\label{appendix:foundation}

\paragraph{Purpose and scope.}
The paper projects use a shared Lean library, called \texttt{Foundation}, for
reusable interfaces and supporting results in convex analysis and
optimization.  The library addresses mathematical components that are only
partially available in \texttt{mathlib} or whose \texttt{mathlib} APIs do not match the notation
used in the source papers.  It is a supporting library rather than a monolithic
import: each paper imports only the modules required by its own definitions and
proofs.

\paragraph{Scale and organization.}
The released version of \texttt{Foundation} contains 61 Lean modules, 48,005 lines of
source code (37,806 non-comment and nonblank lines), and approximately 1,257 top-level
declarations (definitions, lemmas, theorems, structures, classes, instances,
and abbreviations).  The modules are
organized into the 13 thematic namespaces shown in Table~\ref{tab:foundation-themes}.
The largest components are convex analysis, fixed-point theory,
differentiation and ordinary differential equations, and Hilbert space
analysis.

\begin{table}[h]
\centering
\scriptsize
\setlength{\tabcolsep}{3pt}
\begin{tabular}{lrrrp{4.2cm}}
\toprule
Namespace & Modules & Lines & Code lines & Main contents\\
\midrule
Analysis & 6 & 3,277 & 2,574 & Extended-real functions, integrals, sequences, semicontinuity, limiting subdifferentials\\
Convex & 17 & 13,416 & 10,588 & Convexity, Fenchel--Rockafellar duality, subdifferentials, proximal/Moreau operators\\
Differentiation & 4 & 6,699 & 5,816 & Gâteaux and Lipschitz differentiation, Moreau gradients, ODEs\\
FixedPoint & 6 & 6,462 & 4,964 & Nonexpansive, averaged, quasi-nonexpansive and inertial maps\\
Hilbert & 5 & 4,919 & 3,578 & Hilbert geometry, orthogonality, weak topology and convergence\\
Metric & 1 & 577 & 419 & Metric fixed-point results\\
Operator & 3 & 1,323 & 934 & Linear operators, distances and bifunctions\\
Optimization & 3 & 697 & 531 & Argmin, descent structures, and Kurdyka--Lojasiewicz results\\
Order & 1 & 254 & 142 & Extended-real order structures\\
Projection & 1 & 2,316 & 1,657 & Best approximation and projection results\\
Separation & 1 & 1,018 & 766 & Separating hyperplanes\\
SetValued & 10 & 5,562 & 4,823 & Monotone/maximal-monotone operators, enlargements, product normal cones, resolvents, Minty theory\\
Topology & 3 & 1,485 & 1,014 & Topological and sequence-space utilities\\
\midrule
\textbf{Total} & \textbf{61} & \textbf{48,005} & \textbf{37,806} & \\
\bottomrule
\end{tabular}
\caption{Organization of the shared \texttt{Foundation} library. ``Lines'' counts all
source lines, including comments and blank lines; ``Code lines'' excludes comments and
blank lines.}
\label{tab:foundation-themes}
\end{table}

\paragraph{Dependency accounting.}
Each paper project is compiled against the same pinned versions of Lean,
\texttt{mathlib}, and \texttt{Foundation}.  We distinguish three dependency layers:
\texttt{mathlib} provides general mathematical infrastructure, \texttt{Foundation} provides
reusable optimization and convex analysis interfaces, and each paper project contains paper-specific declarations.  The
statistics below count explicit imports of modules whose names begin with
\texttt{Foundation}; they therefore measure library-level reuse.  Declaration-level dependencies can be obtained from Lean's
environment and are kept as part of the released provenance data.

Across the 20 papers, 38 of the 61 \texttt{Foundation} modules are directly imported by
at least one project.  The number of directly imported modules per paper ranges
from 0 to 14 (mean 7.90).  The most widely reused interfaces are
\texttt{Hilbert.WeakConvergence} (14 papers),
\texttt{Analysis.ExtendedRealFunctions} (12),
\texttt{Convex.Subdifferential} (10),
\texttt{SetValued.Monotone} and \texttt{SetValued.Resolvent} (8 each), and
\texttt{FixedPoint.Nonexpansive}, \texttt{Hilbert.WeakTopology}, and
\texttt{SetValued.Basic} (7 each).
Table~\ref{tab:foundation-usage} gives the complete module-level usage summary,
and Table~\ref{tab:paper-foundation} gives the paper-level import summary.

\begin{table}[h]
\centering
\scriptsize
\begin{tabular}{lr@{\quad}lr}
\toprule
Module & Papers & Module & Papers\\
\midrule
\texttt{Hilbert.WeakConvergence} & 14 & \texttt{Analysis.ExtendedRealFunctions} & 12\\
\texttt{Convex.Subdifferential} & 10 & \texttt{SetValued.Monotone} & 8\\
\texttt{SetValued.Resolvent} & 8 & \texttt{FixedPoint.Nonexpansive} & 7\\
\texttt{Hilbert.WeakTopology} & 7 & \texttt{SetValued.Basic} & 7\\
\texttt{Convex.WeakClosed} & 6 & \texttt{Projection.BestApproximation} & 6\\
\texttt{Convex.ProximalOperator} & 5 & \texttt{Differentiation.Gateaux} & 5\\
\texttt{SetValued.Minty} & 5 & \texttt{Convex.StronglyConvex} & 4\\
\texttt{Differentiation.ODE} & 4 & \texttt{Hilbert.Basic} & 4\\
\texttt{Optimization.Argmin} & 4 & \texttt{SetValued.MaximalMonotone} & 4\\
\texttt{Differentiation.MoreauGradient} & 3 & \texttt{Operator.Linear} & 3\\
\texttt{SetValued.Enlargement} & 3 & \texttt{SetValued.MaximalMonotoneSum} & 3\\
\texttt{Analysis.Semicontinuity} & 2 & \texttt{Convex.ArgminExistence} & 2\\
\texttt{Convex.Basic} & 2 & \texttt{Convex.Conjugate} & 2\\
\texttt{Convex.MoreauEnvelope} & 2 & \texttt{Convex.NormalCone} & 2\\
\texttt{FixedPoint.Algorithms} & 2 & \texttt{Metric.FixedPoint} & 2\\
\texttt{SetValued.Yosida} & 2 & \texttt{Topology.Sequences} & 2\\
\texttt{Analysis.LimitingSubdifferential} & 1 & \texttt{Convex.Coercive} & 1\\
\texttt{FixedPoint.Averaged} & 1 & \texttt{Operator.Distance} & 1\\
\texttt{Optimization.KurdykaLojasiewicz} & 1 & \texttt{Separation.Hyperplanes} & 1\\
\bottomrule
\end{tabular}
\caption{Number of \textsc{OptPaper-20} projects that directly import each \texttt{Foundation}
module. Modules with zero direct imports are omitted.}
\label{tab:foundation-usage}
\end{table}

\paragraph{How the papers use the library.}
The dynamical-system papers primarily use Hilbert-space geometry, weak topology,
weak convergence, differentiation, and ODE interfaces.  Proximal and
forward--backward papers use convexity, subdifferentials, proximal operators,
resolvents, and monotone-operator modules.  Fixed-point and projection papers
reuse nonexpansive, averaged, fixed-point, and best-approximation interfaces.
Extended-real functions, semicontinuity, argmin, and coercivity form the common
variational-analysis layer.  Each paper imports a different subset according to
its assumptions and proof structure.

\paragraph{Content boundary.}
\texttt{Foundation} is a shared supporting library, rather than a collection of the
20 dataset papers' results.  It contains generic mathematical interfaces,
including parameterized inertial or fixed-point schemas and general supporting
lemmas, but excludes the paper-specific algorithm update rules, parameter
regimes, and principal convergence theorems of the 20 dataset papers.   Thus, \texttt{Foundation} may provide a general
result about a nonexpansive map or a resolvent, while a concrete FISTA,
inertial forward--backward, proximal, or dynamical-system algorithm and its
main convergence statement are formalized in the relevant paper project.  The
import table measures module-level reuse; it is not by itself a theorem-level
proof that no paper-specific declaration occurs in \texttt{Foundation}.

\begin{table}[h]
\centering
\scriptsize
\begin{tabular}{lrrr}
\toprule
Paper & Lean files & LOC & \texttt{Foundation} modules\\
\midrule
\citet{abbas2015dynamical} & 38 & 10,819 & 10\\
\citet{alvarez2002second} & 33 & 11,752 & 12\\
\citet{alvarez2001inertial} & 8 & 1,428 & 6\\
\citet{Alvarez2000minimizing} & 22 & 7,312 & 9\\
\citet{Alvarez2003weakconverge} & 14 & 4,315 & 9\\
\citet{attouch2014dynamical} & 42 & 6,016 & 9\\
\citet{adly2023accelerated} & 25 & 4,892 & 6\\
\citet{ait2024second} & 41 & 7,301 & 9\\
\citet{attouch2024fast} & 20 & 4,444 & 10\\
\citet{adly2020finite} & 31 & 4,965 & 7\\
\citet{AttouchRedont2007Newclass} & 25 & 8,164 & 4\\
\citet{BajiCabot2006Inertial} & 32 & 15,639 & 14\\
\citet{FISTA2009} & 20 & 2,210 & 0\\
\citet{Bot2016inertialTseng} & 24 & 6,901 & 8\\
\citet{burachik2001robustness} & 19 & 1,684 & 12\\
\citet{dong2017strong} & 14 & 1,943 & 9\\
\citet{Guler1992proximal} & 35 & 2,671 & 4\\
\citet{lorenz2015inertial} & 24 & 4,332 & 9\\
\citet{MoudafiOlny2003converge} & 11 & 2,015 & 7\\
\citet{nesterov1983method} & 16 & 2,107 & 4\\
\bottomrule
\end{tabular}
\caption{Direct \texttt{Foundation} imports in the 20 paper projects. LOC includes
comments and blank lines; the final column counts distinct \texttt{Foundation} modules
imported anywhere under the paper directory.}
\label{tab:paper-foundation}
\end{table}

\paragraph{Validation and versioning.}
\texttt{Foundation} is compiled as an independent Lean library under the same toolchain
and \texttt{mathlib} revision used for the paper projects.  The released snapshot is
identified by its repository revision and is frozen for the reported
experiments.  Every paper project is rebuilt from a clean environment after
integration.  The final audit checks for compilation errors, proof placeholders
(``\texttt{sorry}'' and ``\texttt{admit}''), and unauthorized axioms.  Results
introduced after the frozen snapshot, if any, are reported as paper-local
extensions and are not counted as part of the shared \texttt{Foundation} statistics.

\section{Downstream Theorem Interface Refinement}
\label{appendix:din-downstream}

The reformulation in Figure~\ref{fig:statement-comparison} also simplifies
the downstream convergence theorem. The original theorem carries the
independently supplied gradient, Hessian, velocity, and acceleration maps,
together with hypotheses identifying these maps with the corresponding
derivatives. The revised theorem takes the potential and trajectory directly:
the derivative APIs are used in the trajectory definition, while the required
regularity of the trajectory is stated explicitly. The following listing shows
the two interfaces for the DIN weak-convergence theorem from
\citet{alvarez2002second}.

\begin{figure}[h]
\centering
\begingroup
\definecolor{downBefore}{HTML}{536A82}
\definecolor{downAfter}{HTML}{287B70}
\lstdefinestyle{downstreamLean}{
  basicstyle=\ttfamily\tiny,
  keywordstyle=\color{blue!55!black}\bfseries,
  columns=fullflexible, keepspaces=true, showstringspaces=false,
  breaklines=true, morekeywords={structure,where,theorem,by,Prop},
  literate={ℝ}{{$\mathbb{R}$}}1 {Φ}{{$\Phi$}}1 {α}{{$\alpha$}}1
           {β}{{$\beta$}}1 {∀}{{$\forall$}}1 {∃}{{$\exists$}}1
           {∧}{{$\land$}}1 {≤}{{$\leq$}}1 {→}{{$\rightarrow$}}1
           {∈}{{$\in$}}1 {•}{{$\bullet$}}1
           {‖}{{$\lVert$}}1,
  aboveskip=1pt, belowskip=1pt
}
\tcbset{
  colback=black!1, boxrule=0.6pt, arc=2pt, valign=top,
  equal height group=dinDownstream,
  left=5pt,right=5pt,top=5pt,bottom=5pt,
  fonttitle=\scriptsize\bfseries, coltitle=white,
  before skip=0pt,after skip=0pt
}
\begin{minipage}[t]{0.49\linewidth}\vspace{0pt}
\begin{tcolorbox}[colframe=downBefore,
  title={(a) Before: expanded interface}]
\begin{lstlisting}[style=downstreamLean]
structure AssumptionH (H : Type*)
    [NormedAddCommGroup H] [InnerProductSpace ℝ H]
    [CompleteSpace H] (Φ : H → ℝ)
    (gradΦ : H → H) (hessΦ : H → H →L[ℝ] H) : Prop where
  bounded_below : ∃ m : ℝ, ∀ x, m ≤ Φ x
  twice_continuously_differentiable : ContDiff ℝ 2 Φ
  hessian_is_derivative : ∀ x, HasFDerivAt gradΦ (hessΦ x) x
  hessian_lipschitz_on_bounded :
    ∀ R : ℝ, 0 < R → ∃ L : ℝ, 0 < L → ∀ x y,
      ‖x‖ ≤ R → ‖y‖ ≤ R →
      ‖hessΦ x - hessΦ y‖ ≤ L * ‖x - y‖

theorem din_weakly_converges_to_argmin
    {H : Type*} [NormedAddCommGroup H]
    [InnerProductSpace ℝ H] [CompleteSpace H]
    (Φ : H → ℝ) (gradΦ : H → H)
    (hessΦ : H → H →L[ℝ] H) (α β : ℝ)
    (x v acc : ℝ → H) (hΦ : AssumptionH H Φ gradΦ hessΦ)
    (hgradΦ : ∀ q, HasFDerivAt Φ (innerSL ℝ (gradΦ q)) q)
    (hhessΦ : ∀ R : ℝ, 0 < R → ∃ L : ℝ, 0 < L ∧ ∀ p q,
      ‖p‖ ≤ R → ‖q‖ ≤ R →
      ‖hessΦ p - hessΦ q‖ ≤ L * ‖p - q‖)
    (hconvex : ConvexOn ℝ Set.univ Φ)
    (hmin : {z | ∀ y, Φ z ≤ Φ y}.Nonempty)
    (hsol : IsDINTrajectory H gradΦ hessΦ α β x v acc) :
    (∀ z, (∀ y, Φ z ≤ Φ y) → ∃ d : ℝ,
      Tendsto (fun t => ‖x t - z‖) atTop (nhds d)) ∧
      ∃ z, (∀ y, Φ z ≤ Φ y) ∧
        Tendsto (fun t => toWeakSpace ℝ H (x t)) atTop
          (nhds (toWeakSpace ℝ H z)) := by
  ...
\end{lstlisting}
\end{tcolorbox}
\end{minipage}\hfill
\begin{minipage}[t]{0.49\linewidth}\vspace{0pt}
\begin{tcolorbox}[colframe=downAfter,
  title={(b) After: compact interface}]
\begin{lstlisting}[style=downstreamLean]
structure AssumptionH (H : Type*)
    [NormedAddCommGroup H] [InnerProductSpace ℝ H]
    [CompleteSpace H] (Φ : H → ℝ) : Prop where
  bounded_below : ∃ m : ℝ, ∀ x, m ≤ Φ x
  twice_continuously_differentiable : ContDiff ℝ 2 Φ
  hessian_lipschitz_on_bounded :
    ∀ R : ℝ, 0 < R → ∃ L : ℝ, 0 < L ∧ ∀ x y,
      ‖x‖ ≤ R → ‖y‖ ≤ R →
      ‖fderiv ℝ (gradient Φ) x -
        fderiv ℝ (gradient Φ) y‖ ≤ L * ‖x - y‖

theorem din_weakly_converges_to_argmin
    {H : Type*} [NormedAddCommGroup H]
    [InnerProductSpace ℝ H] [CompleteSpace H]
    (Φ : H → ℝ) (α β : ℝ) (x : ℝ → H)
    (hΦ : AssumptionH H Φ) (hx : ContDiff ℝ 2 x)
    (hconvex : ConvexOn ℝ Set.univ Φ)
    (hmin : (Optimization.argminSet Φ).Nonempty)
    (hsol : IsDINTrajectory H Φ α β x) :
    (∀ z, z ∈ Optimization.argminSet Φ → ∃ d : ℝ,
      Tendsto (fun t => ‖x t - z‖) atTop (nhds d)) ∧
      ∃ z, z ∈ Optimization.argminSet Φ ∧
        Tendsto (fun t => toWeakSpace ℝ H (x t)) atTop
          (nhds (toWeakSpace ℝ H z)) := by
  ...
\end{lstlisting}
\end{tcolorbox}
\end{minipage}
\endgroup
\caption{Downstream interface comparison for the DIN weak convergence
theorem. The revised statement removes repeated map-identification hypotheses,
uses the shared \texttt{argminSet} abstraction, and makes the trajectory
regularity assumption explicit. Proof bodies are omitted.}
\label{fig:din-downstream-comparison}
\end{figure}

The main benefit is interface locality. In the original version, every theorem
that uses the trajectory must thread the auxiliary maps and separately restate
their compatibility with $\Phi$. In the revised version these relationships
are encoded once in the definitions used by the trajectory, while
\texttt{ContDiff} on $x$ records the regularity needed to interpret the second
derivative. The theorem statement is consequently shorter and its hypotheses
are closer to the mathematical objects they constrain. This refactoring does
not remove mathematical obligations: derivative regularity and the
identification of the gradient and Hessian are still required, but they are
represented through shared APIs and reusable assumptions rather than repeated
paper-local parameters.

\section{Statement Revisions During Proof Construction}
\label{appendix:statement-revisions}

\subsection{Case 1: recovering assumptions used in the source proof}
\label{appendix:moudafi-assumptions}

During proof construction for Theorem~2.1 of
\citet{MoudafiOlny2003converge}, the statement was augmented with
\texttt{hS\_nonempty} and \texttt{h\_norm\_diff\_tendsto}
(Figure~\ref{fig:moudafi-assumption-repair}). We regard these additions as
recovering assumptions from the source proof, based on the following evidence.

\paragraph{Existence of a solution.}
The proof begins by fixing a point in $S=(A+B)^{-1}(0)$
(source p.~451). The hypothesis \texttt{hS\_nonempty} makes this implicit
existence assumption explicit and supports the later application of Opial's
lemma; it does not assume convergence of the iterates.

\paragraph{Asymptotic regularity.}
The proof later invokes ``by hypothesis'' the condition
$\|x_k-x_{k-1}\|\to 0$ (source p.~453), although this condition is not listed
in the theorem statement. The hypothesis \texttt{h\_norm\_diff\_tendsto}
expresses the same condition with a shifted index. Our assessment follows
the theorem together with its proof: the revision makes this dependency
explicit, but does not establish that it follows from the theorem's other
listed assumptions. This distinguishes recovery of a source-supported
condition from adding an unsupported restriction merely to complete a proof.

\begin{figure}[htbp]
\centering
\begingroup
\definecolor{assumptionAccent}{HTML}{287B70}
\lstdefinestyle{leanAssumptionRepair}{
  basicstyle=\ttfamily\scriptsize,
  keywordstyle=\color{blue!55!black}\bfseries,
  commentstyle=\color{black!45},
  morekeywords={theorem,fun,by},
  morecomment=[l]{--},
  alsoletter={_},
  emph={hS_nonempty,h_norm_diff_tendsto},
  emphstyle=\color{assumptionAccent}\bfseries,
  literate={ℝ}{{$\mathbb{R}$}}1
           {ℕ}{{$\mathbb{N}$}}1
           {β}{{$\beta$}}1
           {θ}{{$\theta$}}1
           {ε}{{$\varepsilon$}}1
           {δ}{{$\delta$}}1
           {→}{{$\rightarrow$}}1
           {↦}{{$\mapsto$}}1
           {∀}{{$\forall$}}1
           {∃}{{$\exists$}}1
           {∈}{{$\in$}}1
           {∧}{{$\land$}}1
           {≤}{{$\leq$}}1
           {≥}{{$\geq$}}1
           {•}{{$\bullet$}}1
           {‖}{{$\Vert$}}1
           {𝓝}{{$\mathcal{N}$}}1,
  columns=fullflexible,
  keepspaces=true,
  showstringspaces=false,
  breaklines=true,
  aboveskip=0pt,
  belowskip=0pt
}
\begin{tcolorbox}[
  colback=black!1,colframe=black!50,
  boxrule=0.6pt,arc=2pt,
  left=7pt,right=7pt,top=3pt,bottom=3pt,
  fonttitle=\small\bfseries,
  title={Theorem 2.1: statement after assumption recovery}]
\begin{lstlisting}[style=leanAssumptionRepair]
theorem convergence_of_splitting_inertial_proximal
    {A : SetValuedOperator H H} {B : H → H}
    (hA_max : IsMaximalMonotone A)
    (hB_max : IsMaximalMonotone (svEmbed B))
    {β : ℝ} (hB_coco : IsCocoerciveWRT A B β)
    (hB_lipschitz : ∀ x y : H, ‖B x - B y‖ ≤ (1 / β) * ‖x - y‖)
    (x : ℕ → H) (θ lam ε : ℕ → ℝ) (y : ℕ → H)
    (hy : ∀ k, 1 ≤ k → y k = x k + θ k • (x k - x (k - 1)))
    (h_alg : ∀ k, 1 ≤ k →
      (0 : H) ∈ (lam k • enlargement A (ε k) (x (k + 1))) +
        {x (k + 1) - y k} + (lam k • svEmbed B (x k)))
    (h_delta : ∃ (δ : ℝ) (lamUnder : ℝ), 0 < lamUnder ∧ 0 < δ ∧
      ∀ k, 1 ≤ k → lamUnder ≤ lam k ∧ lam k ≤ 2 * β - δ)
    (h_theta : ∃ (θmax : ℝ), 0 ≤ θmax ∧ θmax < 1 ∧
      ∀ k, 1 ≤ k → 0 ≤ θ k ∧ θ k ≤ θmax)
    (h_eps_sum : Summable ε) (h_eps_pos : ∀ k, ε k ≥ 0)
    (h_summable : Summable fun k : ℕ ↦
      θ (k + 1) * ‖x (k + 1) - x k‖ ^ 2)
    (hS_nonempty : ((A + svEmbed B).zeros).Nonempty)
    (h_norm_diff_tendsto :
      Tendsto (fun n : ℕ => ‖x (n+1) - x n‖) atTop (𝓝 0)) :
    ∃ (xstar : H), xstar ∈ (A + svEmbed B).zeros ∧
      Tendsto (fun k ↦ toWeakSpace ℝ H (x k)) atTop
        (𝓝 (toWeakSpace ℝ H xstar)) := by
  -- Proof omitted.
\end{lstlisting}
\end{tcolorbox}
\endgroup
\caption{Revised statement for Theorem~2.1 of
\citet{MoudafiOlny2003converge}. The two hypothesis names highlighted in green
identify additions made during proof construction. The surrounding real
Hilbert-space context and proof body are omitted.}
\label{fig:moudafi-assumption-repair}
\end{figure}

\subsection{Case 2: assuming properties that should be derived}
\label{appendix:dry-friction-assumptions}

The second example concerns Theorem~3.1 of
\citet{adly2023accelerated}, represented by \texttt{thm\_2}. For a convex potential and a solution of (DRYGRAD), the
source establishes that the speed is nonincreasing and
$\|\dot z(t)\|=o(1/t)$. Figure~\ref{fig:dry-friction-assumption-repair}
compares the statement generated by \textsc{Paper2Lean} with the repaired
statement. 
We judge the \textsc{Paper2Lean} output to be unfaithful because it treats intermediate properties as additional premises rather than deriving them from the source assumptions.

In particular, \texttt{hvelocity\_limit} directly assumes $\dot z(t)\to 0$,
and \texttt{hA\_nonneg} supplies the nonnegative gradient-derivative pairing
used in the velocity estimate. Both are consequences of the source
assumptions in the proof: the former is obtained from the integrability result
for the velocity, while the latter follows from convexity and the derivative
of the gradient. The more substantial issue is the auxiliary vector field
\texttt{V}. In the proof, \texttt{V} is constructed rather than passed
as an input, namely
\[
  V(x)=-\operatorname{softThresholder}_{r}(\nabla f(x)).
\]
The proof then derives its local Lipschitz property from the nonexpansiveness
of soft thresholding and the local Lipschitz property of $\nabla f$; proves
\texttt{hV\_zero}, the equivalence
$V(x)=0\ \Longleftrightarrow\ \|\nabla f(x)\|\le r$; and derives
\texttt{hvelocity} from the DRYGRAD subdifferential inclusion. Finally,
\texttt{hstationary} is proved using uniqueness for the ODE with this locally
Lipschitz vector field, rather than assumed as an independent property.
Thus \texttt{V}, its compatibility with the trajectory, its zero-set
characterization, and the stationarity consequence are all intermediate
objects or lemmas obtained inside the proof. Treating them as hypotheses
changes proof obligations into extra restrictions on the theorem and
can make the formalized result strictly weaker than the source theorem.

The repaired development uses \texttt{gradient f} and \texttt{deriv z}
directly, derives the nonnegative pairing from convexity, and obtains velocity
integrability from the preceding formalized Theorem~2.1. The preceding
development also derives velocity convergence instead of requiring
\texttt{hvelocity\_limit} as an input. This repairs the dependency structure
as well as simplifying the interface. The repaired statement explicitly
assumes $C^2$ regularity of $f$ and $z$; the comparison demonstrates removal
of the additional intermediate premises under that regularity, not
unconditional equivalence of the two displayed statements.

\begin{figure}[h]
\centering
\begingroup
\definecolor{failedAssumption}{HTML}{9A433A}
\definecolor{repairedStatement}{HTML}{287B70}
\lstdefinestyle{leanDryFrictionRepair}{
  basicstyle=\ttfamily\scriptsize,
  keywordstyle=\color{blue!55!black}\bfseries,
  commentstyle=\color{black!45},
  morekeywords={theorem,Type,fun,by},
  morecomment=[l]{--},
  alsoletter={_},
  literate={ℝ}{{$\mathbb{R}$}}1
           {₀}{{$_0$}}1
           {→}{{$\rightarrow$}}1
           {↔}{{$\leftrightarrow$}}1
           {∀}{{$\forall$}}1
           {∃}{{$\exists$}}1
           {∈}{{$\in$}}1
           {∧}{{$\land$}}1
           {≤}{{$\leq$}}1
           {‖}{{$\Vert$}}1,
  columns=fullflexible,
  keepspaces=true,
  showstringspaces=false,
  breaklines=true,
  aboveskip=0pt,
  belowskip=0pt
}
\tcbset{
  colback=black!1,boxrule=0.6pt,arc=2pt,
  left=7pt,right=7pt,top=3pt,bottom=3pt,
  fonttitle=\small\bfseries,coltitle=white,
  before skip=3pt,after skip=3pt
}
\begin{tcolorbox}[colframe=failedAssumption,
  title={(a) Automatic statement: additional intermediate premises}]
\begin{lstlisting}[style=leanDryFrictionRepair,
  emph={hvelocity_limit,hA_nonneg,hstationary,hV_zero,hvelocity},
  emphstyle=\color{failedAssumption}\bfseries]
theorem theorem_3_1 {H : Type*} [NormedAddCommGroup H] [InnerProductSpace ℝ H]
    [CompleteSpace H]
    {t₀ r : ℝ} (hr : 0 < r) (f : H → ℝ) (fgrad : H → H) (z dz : ℝ → H)
    (hconvex : ConvexOn ℝ Set.univ f)
    (hgradient : ∀ x v, fderiv ℝ f x v = inner ℝ (fgrad x) v)
    (hsol : IsDrygradTrajectory r hr fgrad z dz t₀)
    (hf_lower : ∃ m : ℝ, ∀ x : H, m ≤ f x) (hfgrad_lip : LocallyLipschitz fgrad)
    (hinitial : r < ‖fgrad (z t₀)‖)
    (hchain : ∀ t ∈ Set.Ici t₀,
      HasDerivAt (fun s => f (z s)) (inner ℝ (fgrad (z t)) (dz t)) t)
    (hcomp_cont : ContinuousOn (fun t => f (z t)) (Set.Ici t₀))
    (V : H → H) (hV_lip : LocallyLipschitz V)
    (hV_zero : ∀ x, V x = 0 ↔ ‖fgrad x‖ ≤ r)
    (hvelocity : ∀ t ∈ Set.Ici t₀, dz t = V (z t))
    (hvelocity_limit : Filter.Tendsto dz Filter.atTop (nhds 0))
    (hstationary : ∀ S ∈ Set.Ici t₀, V (z S) = 0 →
      ∀ t ∈ Set.Ici S, z t = z S)
    (ddz : ℝ → H) (A : ℝ → H → H)
    (hdz_deriv : ∀ t ∈ Set.Ici t₀, HasDerivAt dz (ddz t) t)
    (hgrad_deriv : ∀ t ∈ Set.Ici t₀,
      HasDerivAt (fun s => fgrad (z s)) (A t (dz t)) t)
    (hA_nonneg : ∀ t ∈ Set.Ici t₀, 0 ≤ inner ℝ (A t (dz t)) (dz t)) :
    AntitoneOn (fun t => ‖dz t‖) (Set.Ici t₀) ∧
      Filter.Tendsto (fun t : ℝ => t * ‖dz t‖) Filter.atTop (nhds 0) := by
  -- Proof omitted.
\end{lstlisting}
\end{tcolorbox}
\begin{tcolorbox}[colframe=repairedStatement,
  title={(b) Human-repaired statement}]
\begin{lstlisting}[style=leanDryFrictionRepair]
theorem theorem_3_1 {H : Type*} [NormedAddCommGroup H] [InnerProductSpace ℝ H]
    [CompleteSpace H]
    {t₀ r : ℝ} (hr : 0 < r) (f : H → ℝ) (z : ℝ → H)
    (hf : ContDiff ℝ 2 f)
    (hconvex : ConvexOn ℝ Set.univ f)
    (hsol : IsDrygradTrajectory r f z t₀)
    (hf_lower : ∃ m : ℝ, ∀ x : H, m ≤ f x)
    (hgradient_lip : LocallyLipschitz (gradient f))
    (hz_contdiff : ContDiff ℝ 2 z)
    (hinitial : r < ‖gradient f (z t₀)‖) :
    AntitoneOn (fun t => ‖deriv z t‖) (Set.Ici t₀) ∧
      Filter.Tendsto (fun t : ℝ => t * ‖deriv z t‖) Filter.atTop (nhds 0) := by
  -- Proof omitted.
\end{lstlisting}
\end{tcolorbox}
\endgroup
\caption{Statement repair for Theorem~3.1 of
\citet{adly2023accelerated}. The highlighted names identify
selected intermediate premises and the convexity assumption in the
automatic version. Each panel uses its corresponding version of
\texttt{IsDrygradTrajectory}; proof bodies are omitted.}
\label{fig:dry-friction-assumption-repair}
\end{figure}

\section{Limitations}

First, the quality and efficiency of downstream formalization depend strongly on the extracted natural-language items, so the extraction stage remains a limitation. Missing assumptions scattered throughout a paper can cause proof construction to stall and waste tokens. Second, we have not compared our approach with other pipelines or studied how different harness architectures affect performance. Future work will evaluate a broader range of pipelines and harness designs. For example, \texttt{DeepSeek-V4.1-Flash} may perform better with the DeepSeek Harness~\citep{deepseek-harness2026}.

\end{document}